%%%%%%%%%%%%%%%%%%%%%%%%%%%%%%%%%%%%%%%%%%%%%%%%%%%%%%%%%%
%%%%%%%%%%%%%%%%%%%%%%%%%%%%%%%%%%%%%%%%%%%%%%%%%%%%%%%%%%
%%%    This is the AMS-LaTeX file:
%%
%%     Colli-Gilardi-Signori-Sprekels CGSS8
%%    sixth-order Cahn--Hilliard + Brinkman
%%%%%%%%%%%%%%%%%%%%%%%%%%%%%%%%%%%%%%%%%%%%%%%%%%%%%%%%%%
%%%%%%%%%%%%%%%%%%%%%%%%%%%%%%%%%%%%%%%%%%%%%%%%%%%%%%%%%%

\def\LaTeX{%
  \let\Begin\begin
  \let\End\end
  \let\salta\relax
  \let\finqui\relax
  \let\futuro\relax}

\def\UK{\def\our{our}\let\sz s}
\def\USA{\def\our{or}\let\sz z}

\UK 
%\USA

%%%%%%%%%%%%%%%%%%%%%%%%%%%%%%%%%

% scegliere fra \TeX e \LaTeX  e fra  \UK oppure \USA

%\TeX
\LaTeX

%\UK
\USA

%%%%%%%%%%%%%%%%%%%%%%%%%%%%%%%%%
%% page layout
%%%%%%%%%%%%%%%%%%%%%%%%%%%%%%%%%

\salta

\documentclass[twoside,12pt]{article}
\setlength{\textheight}{24cm}
\setlength{\textwidth}{16cm}
\setlength{\oddsidemargin}{2mm}
\setlength{\evensidemargin}{2mm}
\setlength{\topmargin}{-15mm}
\parskip2mm

%%%%%%%%%%%%%%%%%%%%%%%%%%%%%%%%%
%% packages
%%%%%%%%%%%%%%%%%%%%%%%%%%%%%%%%%

%\usepackage{color}
\usepackage[usenames,dvipsnames]{color}
\usepackage{amsmath}
\usepackage{amsthm}
\usepackage{amssymb,bbm}
\usepackage[mathcal]{euscript}

\usepackage{cite}
\usepackage{hyperref}
\usepackage[shortlabels]{enumitem}

\usepackage[ulem=normalem,draft]{changes}
%
%		COLORS FOR CORRECTIONS
%
% do the same, please (i.e., don't use the standard {\color{red} text} or similar): 
% just choose the color you prefer in \def\yourname

% EXAMPLE OF USE:  \fredi{I want this to become blue}
%
%IF YOU LATER WANT TO LET THE COLOR DISAPPEAR, ACTIVATE \def\fredi #1{{#1}} BELOW
 
\definecolor{viola}{rgb}{0.3,0,0.7}
\definecolor{ciclamino}{rgb}{0.5,0,0.5}
\definecolor{blu}{rgb}{0,0,0.7}
\definecolor{rosso}{rgb}{0.85,0,0}

%%%%%%%%%%%%%     ADDED BY PIER

%\beamertemplateballitem

\definecolor{dgreen}{rgb}{0,0.6,0.1}
\definecolor{dred}{rgb}{0.8,0,0}
\definecolor{dblue}{rgb}{0,0,0.8}
\definecolor{dbrown}{rgb}{0.5,0.25,0}
\definecolor{dorange}{rgb}{1,0.4,0}

\definecolor{lightgreen}{rgb}{0.8,1,0.8}
\definecolor{newyel}{rgb}{1,1,0.7}
\definecolor{lightcyan}{rgb}{0.8,1,1}
\definecolor{lightpink}{rgb}{1,0.9,0.9}
\definecolor{lightgrey}{rgb}{0.9,0.9,0.9}
\definecolor{lightbrown}{rgb}{0.8,0.8,0.7}

\definecolor{gray}{rgb}{0.5,0.5,0.7}
\definecolor{ddgreen}{rgb}{0,0.4,0.4}
\definecolor{miogreen}{rgb}{0,0.5,0.2}
\definecolor{magenta-mio}{rgb}{1,0,1}
\definecolor{dcyan}{rgb}{0,0.4,0.8}
\definecolor{dmagenta}{rgb}{0.8,0,0.8}
\definecolor{dyellow}{rgb}{0.8,0.8,0}

%%%%%%%%%%%%%    

\def\juerg #1{{\color{black}#1}}
\def\junew #1{{\color{green}#1}}
\def\anold #1{{\color{black}#1}}
\def\an #1{{\color{dblue}#1}}
\def\last #1{{\color{magenta}#1}}
\def\pier #1{{\color{rosso}#1}} 

\def\julast #1{{\color{magenta}#1}}
\def\pcol #1{{\color{rosso}#1}}
\def\last  #1{{\color{rosso}#1}}

\def\juerg #1{#1}
\def\junew #1{#1} 
\def\an #1{#1}
\def\pier #1{#1}

\def\julast #1{#1}
\def\pcol #1{#1}
\def\last #1{#1}

%%%%%%%%%%%%%%%%%%%%%%%%%%%%%%%%%
%% you may adjust the baseline
%%%%%%%%%%%%%%%%%%%%%%%%%%%%%%%%%

%\renewcommand{\baselinestretch}{0.975}

%%%%%%%%%%%%%%%%%%%%%%%%%%%%%%%%%
%% bibliographystyle
%%%%%%%%%%%%%%%%%%%%%%%%%%%%%%%%%

\bibliographystyle{plain}

%%%%%%%%%%%%%%%%%%%%%%%%%%%%%%%%%
%% environments
%%%%%%%%%%%%%%%%%%%%%%%%%%%%%%%%%

%
\newtheorem{theorem}{Theorem}[section]

\finqui

\def\Beq{\Begin{equation}}
\def\Eeq{\End{equation}}

\def\Bthm{\Begin{theorem}}
\def\Ethm{\End{theorem}}

\def\Brem{\Begin{remark}\rm}
\def\Erem{\End{remark}}

\def\Bcenter{\Begin{center}}
\def\Ecenter{\End{center}}
\let\non\nonumber

%%%%%%%%%%%%%%%%%%%%%%%%%%%%%%%%%
%% macros
%%%%%%%%%%%%%%%%%%%%%%%%%%%%%%%%%

% macro salvate

% sottosezioni non numerate

\def\step #1 \par{\medskip\noindent{\bf #1.}\quad}
\def\jstep #1: \par {\vspace{2mm}\noindent\underline{\sc #1 :}\par\nobreak\vspace{1mm}\noindent}

\def\aand{\quad\hbox{and}\quad}
\def\Lip{Lip\-schitz}
\def\Holder{H\"older}

\def\lhs{left-hand side}
\def\rhs{right-hand side}

% versioni inglesi (UK) o americane (USA)

%\def\analyz {analy\sz}

% bold, cal, grass e mathop

\def\multibold #1{\def\arg{#1}%
  \ifx\arg\pto \let\next\relax
  \else
  \def\next{\expandafter
    \def\csname #1#1\endcsname{{\boldsymbol #1}}%
    \multibold}%
  \fi \next}

\def\pto{.}

\def\multical #1{\def\arg{#1}%
  \ifx\arg\pto \let\next\relax
  \else
  \def\next{\expandafter
    \def\csname cal#1\endcsname{{\cal #1}}%
    \multical}%
  \fi \next}

\def\multigrass #1{\def\arg{#1}%
  \ifx\arg\pto \let\next\relax
  \else
  \def\next{\expandafter
    \def\csname grass#1\endcsname{{\mathbb #1}}%
    \multigrass}%
  \fi \next}

% operatori

\def\multimathop #1 {\def\arg{#1}%
  \ifx\arg\pto \let\next\relax
  \else
  \def\next{\expandafter
    \def\csname #1\endcsname{\mathop{\rm #1}\nolimits}%
    \multimathop}%
  \fi \next}

\multibold
qweryuiopasdfghjklzxcvbnmQWERTYUIOPASDFGHJKLZXCVBNM.  % esclusa t per non cambiare \tt

\multical
QWERTYUIOPASDFGHJKLZXCVBNM.

\multigrass
QWERTYUIOPASDFGHJKLZXCVBNM.

\multimathop
diag dist div dom mean meas sign supp .

\def\Span{\mathop{\rm span}\nolimits}

% accorpamenti di formule citate:
% uso  \accorpa {prima}{seconda}
%      \Accorpa\cs prima seconda (con il comodo blank anche dopo)
% NB: \Accorpa definisce \cs come l'accorpamento delle due citazioni
% e scrive sul file.log

\def\accorpa #1#2{\eqref{#1}--\eqref{#2}}
\def\Accorpa #1#2 #3 {\gdef #1{\eqref{#2}--\eqref{#3}}%
  \wlog{}\wlog{\string #1 -> #2 - #3}\wlog{}}

% macro comode

\def\separa{\noalign{\allowbreak}}

\def\somma #1#2#3{\sum_{#1=#2}^{#3}}

\def\graffe #1{\mathopen\{#1\mathclose\}}
\def\<#1>{\mathopen\langle #1\mathclose\rangle}
\def\norma #1{\mathopen \| #1\mathclose \|}

\def\aeQ{\checkmmode{a.e.\ in~$Q$}}

\def\aet{\checkmmode{a.e.\ in~$(0,T)$}}
\def\aat{\checkmmode{for a.a.\ $t\in(0,T)$}}

\def\cpto{\,\cdot\,}

\def\iot {\int_0^t}

\def\intQt{\int_{Q_t}}
\def\intQ{\int_Q}
\def\iO{\int_\Omega}

\def\dt{\partial_t}
\def\dn{\partial_{\nn}}
\def\ddt{\frac d{dt}}

\let\emb\hookrightarrow
\def\cpto{\,\cdot\,}

\def\checkmmode #1{\relax\ifmmode\hbox{#1}\else{#1}\fi}

% insiemi numerici

\let\erre\grassR
\let\enne\grassN
\def\erren{\erre^n}

% spazi di funzioni a valori vettoriali su [0,T], [0,t], [0,s], [0,+\infty), [\delta,T]

% Come ricordare: in generale i simboli L H W  C da soli per gli spazi su (0,T)
% gli stessi raddoppiati per (0,+\infty)
% aggiunta di t o s al simbolo per (0,t) e (0,s)
% aggiunta di d al simbolo semplice o doppio per intervalli (\delta,T) e (\delta,+\infty)
% il simbolo C e i suoi derivati mettono le quadre anziche' le tonde

% Esempi   \L2V   \L\infty\Vp   \W{1,1}H   \C0H   \LL2V   \calC0\Vp   \Ld2V  \calCdH

\def\genspazio #1#2#3#4#5{#1^{#2}(#5,#4;#3)}
\def\spazio #1#2#3{\genspazio {#1}{#2}{#3}T0}

\def\L {\spazio L}
\def\H {\spazio H}

\def\C #1#2{C^{#1}([0,T];#2)}

% spazi di funzioni su \Omega, \Gamma, Q e \Sigma

\def\Lx #1{L^{#1}(\Omega)}
\def\Hx #1{H^{#1}(\Omega)}
\def\Wx #1{W^{#1}(\Omega)}

\def\LQ #1{L^{#1}(Q)}

\def\CS #1{C^{#1}(\Sigma)}

\def\Luno{\Lx 1}
\def\Ldue{\Lx 2}
\def\Linfty{\Lx\infty}

\def\Huno{\Hx 1}
\def\Hdue{\Hx 2}

% simboli in bold

\def\LLx #1{\LL^{#1}(\Omega)}
\def\HHx #1{\HH^{#1}(\Omega)}

% lettere greche

%\let\badtheta\theta
%\let\theta\vartheta
%\let\badeps\epsilon
\let\eps\varepsilon
\let\badphi\phi
\let\phi\varphi
\def\zz{\boldsymbol\zeta}

\let\TeXchi\chi                         % new \chi, exactly on the baseline
\newbox\chibox
\setbox0 \hbox{\mathsurround0pt $\TeXchi$}
\setbox\chibox \hbox{\raise\dp0 \box 0 }
\def\chi{\copy\chibox}

% quadratino di fine dimostrazione

% abbreviazioni specifiche del lavoro

\def\zbar{\overline z}
\def\mubar{\overline\mu}

\def\ej{e_j}
\def\ei{e_i}
\def\lambdaj{\lambda_j}
\def\Vn{V_n}
\def\phin{\phi_n}
\def\mun{\mu_n}
\def\wn{w_n}
\def\vnj{v_{nj}}
\def\phinj{\phi_{nj}}
\def\munj{\mu_{nj}}
\def\vni{v_{ni}}
\def\phini{\phi_{ni}}
\def\muni{\mu_{ni}}

\def\zn{z_n}

\def\eej{\ee_{0,j}}
\def\eei{\ee_{0,i}}
\def\lambdazj{\lambda_{0,j}}
\def\vvn{\vv_n}
\def\VVz{\VV_{\!0}}
\def\HHz{\HH_0}
\def\VVzn{\VV_{0,n}}
\def\VVzm{\VV_{0,\anold{N}}}
\def\0{\boldsymbol 0}

\let\T\grassT
\let\I\grassI

\def\Pn{\grassP_n}

\def\phiz{\phi_0}

\def\Vp{{V^*}}

\def\soluz{(\vv,\phi,\mu,w)}
\def\soluzn{(\vvn,\phin,\mun)}

\def\normaV #1{\norma{#1}_V}
\def\normaW #1{\norma{#1}_W}

\def\etamin{\eta_*}
\def\etamax{\eta^*}
\def\lambdamin{\lambda_*}
\def\lambdamax{\lambda^*}
\def\mmin{m_*}
\def\mmax{m^*}

\def\CO{C_\Omega}
\def\CK{\calC_K}
\def\CP{\calC_P}
\def\CS{\calC_S}
\def\CE{\calC_E}
\def\Cdelta{\calC_\delta}
\def\cdelta{c_\delta}

\def\cM{c_M}

\def\cdeltaM{c_{\delta,M}}

%%%%%%%%%%%%%%%%%%%%%%%%%%%%%%%%%%%%%%%%%%%%%%%%%%

\usepackage{amsmath}
\DeclareFontFamily{U}{mathc}{}
\DeclareFontShape{U}{mathc}{m}{it}%
{<->s*[1.03] mathc10}{}

\DeclareMathAlphabet{\mathscr}{U}{mathc}{m}{it}

%%%%%%%%%%%%%%%%%%%%%%%%%%%%%%
\Begin{document}
%%%%%%%%%%%%%%%%%%%%%%%%%%%%%%%%%

%%%%%%%%%%%%%%%%%%%%%%%%%%%%%%%%%
%% front page
%%%%%%%%%%%%%%%%%%%%%%%%%%%%%%%%%

%
\title{
On Brinkman flows with curvature-induced phase separation in binary mixtures
}

\author{}
\date{}
\maketitle
\Bcenter
\vskip-1.5cm
{\large\sc Pierluigi Colli$^{(1)}$}\\
{\normalsize e-mail: {\tt pierluigi.colli@unipv.it}}\\
[0.25cm]
{\large\sc Gianni Gilardi$^{(2)}$}\\
{\normalsize e-mail: {\tt gianni.gilardi@unipv.it}}\\
[0.25cm]
{\large\sc Andrea Signori$^{(3)}$}\\
{\normalsize e-mail: {\tt andrea.signori@polimi.it}}\\
[0.25cm]
{\large\sc J\"urgen Sprekels$^{(4)}$}\\
{\normalsize e-mail: {\tt juergen.sprekels@wias-berlin.de}}\\
[.5cm]
$^{(1)}$
{\small Dipartimento di Matematica ``F. Casorati'', Universit\`a di Pavia}\\
{\small and Research Associate at the IMATI -- C.N.R. Pavia}\\
{\small via Ferrata 5, I-27100 Pavia, Italy}\\
[0.25cm]
$^{(2)}$
{\small Dipartimento di Matematica ``F. Casorati'', Universit\`a di Pavia}\\
{\small via Ferrata 5, I-27100 Pavia, Italy}\\
{\small and Istituto Lombardo, Accademia di Scienze e Lettere}\\
{\small via Borgonuovo 25, I-20121 Milano, Italy}
\\[0.25cm]
$^{(3)}$
{\small Dipartimento di Matematica, Politecnico di Milano}\\
{\small via E. Bonardi 9, I-20133 Milano, Italy}
\\
\anold{
{Alexander von Humboldt Research Fellow}}
\\[.25cm] 
$^{(4)}$
{\small Department of Mathematics}\\
{\small Humboldt-Universit\"at zu Berlin}\\
{\small Unter den Linden 6, D-10099 Berlin, Germany}\\
{\small and Weierstrass Institute for Applied Analysis and Stochastics}\\
{\small \junew{Anton-Wilhelm-Amo-Strasse} 39, D-10117 Berlin, Germany}\\[10mm]
\date{}
\Ecenter
\Begin{abstract}
\noindent
The mathematical analysis of diffuse-interface models for multiphase flows has attracted significant attention due to their ability to capture complex interfacial dynamics, including curvature effects, within a unified, energetically consistent framework. In this work, we study a novel Brinkman--Cahn--Hilliard system, coupling a sixth-order phase-field evolution with a Brinkman-type momentum equation featuring variable shear viscosity. The Cahn--Hilliard equation includes a nonconservative source term accounting for mass exchange, and the velocity equation contains a non divergence-free forcing term. We establish the existence of weak solutions in a divergence-free variational framework, and, in the case of constant mobility and shear viscosity, prove uniqueness and continuous dependence on the forcing. Additionally, we analyze the Darcy limit, providing existence results for the corresponding reduced system.


%    35K35 — Higher-order parabolic equations
%
%    35Q35 — PDEs in fluid mechanics
%
%    76D07 — Stokes and related (viscous) flows
%
%    76T99 — Multiphase and complex fluids (miscellaneous)

\End{abstract}

\pagestyle{myheadings}
\newcommand\testopari{\sc Colli -- Gilardi -- Signori -- Sprekels}
\newcommand\testodispari{\sc Brinkman--Cahn--Hilliard model with curvature effects}
\markboth{\testopari}{\testodispari}
%
%%%%%%%%%%%%%%%%%%%%%%%%%%%%%%%%%
%% very beginning
%%%%%%%%%%%%%%%%%%%%%%%%%%%%%%%%%

\section{Introduction}
\label{INTRO}
\setcounter{equation}{0}

Multiphase flow research has been profoundly advanced by diffuse-interface models, which provide a rigorous energetic framework to capture complex interfacial phenomena, including curvature effects. In particular, models based on the Cahn–Hilliard theory and its higher-order variants naturally describe binary mixtures, amphiphilic membranes, and phase separation in soft matter. Such models are especially relevant in biology, where lipid bilayers and cellular membranes undergo curvature-driven transformations, fusion, and budding, processes that are crucial for intracellular organization. When coupled with fluid dynamics, they provide a powerful tool to investigate the interplay between interfacial evolution, transport, and hydrodynamic effects.

Along these lines,  \julast{we investigate in this paper a novel system of equations of Brinkman--Cahn–Hilliard type.} 
Let $\Omega\subset\erre^3$ be a bounded domain with smooth boundary~$\Gamma$.
If $T>0$ is a fixed final time, we set
$Q:=\Omega\times(0,T)$ and $\Sigma:=\Gamma\times(0,T)$ and
\Beq
  \non
  Q_t := \Omega \times (0,t)
  \quad \hbox{for $t\in(0,T]$} \,.
%  \label{defQt}
\Eeq
Then, the differential system under investigation reads
\begin{align}
  & - \div \T(\phi,\vv,p)
  + \lambda(\phi) \vv
  = \mu \nabla\phi
  + \uu
  \aand 
  \div\vv = 0
  & \quad \hbox{in $Q$}\,,
  \label{Iprima}
  \\
  & \dt\phi
  + \vv \cdot \nabla\phi
  - \div(m(\phi)\nabla\mu)
  = S(\phi) 
%  \pier{{}- \sigma \phi}
  & \quad \hbox{in $Q$}\,,
  \label{Iseconda}
  \\
  & - \eps \Delta w
  + \tfrac 1 \eps f'(\phi) w
  + \nu w 
  = \mu
  & \quad \hbox{in $Q$}\,,
  \label{Iterza}
  \\
  & - \eps\Delta\phi
  + \tfrac 1 \eps f(\phi) 
  = w
  & \quad \hbox{in $Q$}\,,
  \label{Iquarta}  
\end{align}
complemented with the boundary and initial conditions
\begin{align}
  & \T(\phi,\vv,p) \nn = \0 
  \aand
  \dn\mu = \dn w = \dn\phi = 0
  & \quad \hbox{on $\Sigma$}\,,
  \label{Ibc}
  \\
  & \phi(0) = \phiz
  & \quad \hbox{in $\Omega$}\,,
  \label{Icauchy}
\end{align}
\Accorpa\Ipbl Iprima Icauchy
where $\nn$ and $\dn$ denote the outward unit normal vector to $\Gamma$ and the corresponding outward normal derivative, respectively, and $\phiz$ is a prescribed function, which serves as the initial datum.

The overall model presents several distinctive features. The velocity field $\vv$ satisfies a Brinkman-type momentum equation \eqref{Iprima}, where the stress tensor $\T$ involves the symmetrized gradient $D\vv$ of $\vv$:
\begin{align*}
	\T(\phi,\vv,p)
	& = \eta(\phi) D \vv
	- p \I \,,
%	\label{defT}
\end{align*}
with phase-dependent positive viscosity $\eta$, $p$ the pressure, $\I\in\erre^{3\times 3}$ the identity matrix, and
\Beq
  D\vv := \frac 12 \, \bigl( \nabla\vv + (\nabla\vv)^{\top} \bigr)
  \label{symgrad}
\Eeq
denoting the symmetrized gradient. This formulation interpolates between the Stokes (formally $\lambda\equiv 0$) and Darcy regimes (formally $\eta\equiv 0$) and is particularly suitable for flows in porous or heterogeneous media.  
\pcol{Since in this work we also rigorously analyze the Darcy limit of the system, obtained by letting the shear viscosity vanish, it is worth emphasizing that the Darcy regime plays a particularly relevant role in applications to biological and porous media, as it effectively captures the flow through highly permeable tissues and fibrous structures where viscous effects are negligible.}
The order parameter $\phi$ represents the local proportion of one of the two components in the binary material and serves
as an order parameter. To simplify the analysis, it is usually normalized in such a way
that the pure states correspond to $\phi=-1$ and  $\phi=1$, while $\{  -1 < \phi  < 1\}$  represents the diffuse
interface layer that occurs in an tubular neighborhood of the interface, with thickness
parameter $\eps  > 0. $
 Its evolution is governed by a \julast{sixth}-order variant of the Cahn--Hilliard-type equation with a source term \julast{of the form}
\[
	S(\phi)= - \sigma \phi + h(\phi),
\]
where $h$ is smooth and bounded, and $\sigma\in\erre$ is a real coefficient. This term accounts for mass exchange and explicitly violates the standard mass-conservation property of $\phi$.  
The chemical potentials $\mu$ and $w$ appearing in equations \eqref{Iterza} and \eqref{Iquarta} correspond to the first variations of the total free energy $\calE$ and the Ginzburg--Landau free energy $\calG$, respectively:
\[
	\mu := \frac{\delta \calE}{\delta \phi}, \qquad w := \frac{\delta \calG}{\delta \phi},
\]
where
\begin{align}
 \calE(\phi) 
 & := 
 \calF(\phi) + \nu \, \calG(\phi)
 =
 \frac 12 \iO  (-\eps \Delta\phi + \tfrac 1 \eps f(\phi))^2
 + \nu  \iO  \Bigl( \frac\eps2 \, |\nabla\phi|^2 + \tfrac 1 \eps F(\phi) \Bigl) \,.
  \label{defE} 
\end{align}
Here, $F$ is an everywhere defined double-well potential, $f$ its derivative, and $\nu$ a \juerg{(not necessarily positive) real constant}. A prototype is the \emph{standard regular potential} 
\Beq
  F_{\rm reg}(s) := \frac 14 \, (s^2-1)^2\,, \quad s\in\erre\,.
  \label{Freg}
\Eeq
Finally, $\uu$ in \eqref{Iprima} represents a forcing term, which can be interpreted as a control variable in a related control problem \pcol{(cf. \cite{CGSS-SIAM, SigW})}, while the nonlinear term $\mu \nabla \phi$ corresponds to the Korteweg force, accounting for capillarity effects in the mixture. 
It is worth mentioning that, due to the regularity of $f$, the boundary condition $\dn w = 0$ in \eqref{Ibc} is equivalent to $\dn \Delta \phi = 0$ on $\Sigma$.  
Overall, the system naturally accommodates a wide range of physical regimes depending on the sign and magnitude of $\nu$ and combines a \julast{sixth}-order Cahn--Hilliard type equation with a Brinkman flow through transport and capillarity effects, the nonconservative source term $S(\phi)$, and the forcing $\uu$.

Let us now frame the model within the literature. A key application concerns curvature-driven phenomena, notably the evolution of amphiphilic bilayer membranes. The classical Helfrich model~\cite{Canham1970, Helfrich1973} describes the bending energy of a membrane $\Gamma_0$ as
\begin{align*}
  {\cal E}_{\rm elastic} = \frac{k}{2} \int_{\Gamma_0} (H - H_0)^2 \, \mathrm{d}S,
\end{align*}
where $H$ is the mean curvature, $H_0$ the spontaneous curvature (often zero), and $k$ the bending rigidity. In a diffuse-interface framework, a phase variable $\phi$ distinguishes interior ($\phi=1$) and exterior ($\phi=-1$) regions, and the Helfrich energy can be approximated by a modified Willmore functional
\begin{align}\label{Willmore}
  {\cal E}_\varepsilon(\phi) = \frac{k}{2\varepsilon} \int_{\Omega} \Bigl(-\varepsilon \Delta \phi + \frac{1}{\varepsilon} (\phi^2 - 1)\phi \Bigr)^2 ,
\end{align}
with $\eps>0$  the interfacial thickness~\cite{DuLiuRyhamWang2005, DuLiuRyhamWang2009}, converging to the sharp-interface limit as $\eps \to 0$. \pcol{Since sharp-interface asymptotics are not considered in this work}, we set $\eps=1$ from now onward for convenience.

It is worth noting that the total energy ${\cal E}$ in our model can be viewed as a higher-order extension of the classical Ginzburg--Landau free energy ${\cal G}$. When $\nu = 0$, ${\cal E}$ reduces to the Willmore functional \eqref{Willmore}, recovering the Canham--Helfrich bending energy of biomembranes. For $\nu > 0$, ${\cal E}$ provides a regularized version of ${\cal G}$, penalizing high curvature and favoring smoother interfaces. This type of regularization has been used, for example, in \cite{TLVW, BCMS} to study anisotropy effects in thin film growth and coarsening processes.  
Conversely, for $\nu < 0$, ${\cal E}$ corresponds to the so-called functionalized Cahn--Hilliard energy, originally introduced for amphiphilic mixtures and later applied to polymers and bilayer membranes, see \cite{GS}. These different regimes illustrate the versatility of the model, capable of capturing a wide range of interfacial phenomena, from membrane elasticity to pattern formation in soft matter.  
Without claiming exhaustiveness, we refer to~\cite{SW, ChengWangWiseYuan2020, ClimentEzquerraGuillenGonzalez2019, ColliLaurencot2011,ColliLaurencot2012, DuLiLiu2007, WuXu2013} for analytical studies and to~\cite{AlandEgererLowengrubVoigt2014, CampeloHernandezMachado2006, Canham1970, DuLiuRyhamWang2005, DuLiuRyhamWang2005b, DuLiuRyhamWang2009, DuLiuWang2004,  LowengrubRatzVoigt2009} for numerical investigations.
The analysis of sixth-order Cahn--Hilliard equations has attracted considerable attention, with applications ranging from the dynamics of oil-water-surfactant mixtures~\cite{PZ1, PZ2, SP}, and the phase-field-crystal model~\cite{GW1, GW2, M1, M2, WW}. 
In the context of applications to optimal control, see also~\cite{SigW, CGSS-SIAM}. Finally, we mention some works on coupling the Cahn–Hilliard equation with Brinkman-type flow~\cite{E_thesis, EG_asy, EL, KS2, CKSS, BCG, ContiG}.

\pcol{The main goal of this work is to establish the well-posedness of the initial-boundary value problem \Ipbl. This is achieved by rewriting the system, especially the velocity equation, in a variational framework where the pressure does not appear explicitly.}
The structure of the paper is as follows. In Section~\ref{STATEMENT}, we introduce the notation, state the precise assumptions, and formulate the problem under consideration. Section~\ref{EXISTENCE} is devoted to the proof of the existence theorem. 
In Section~\ref{PROOFDARCY},
\pcol{we establish the Darcy limit by letting the shear viscosity tend to zero, thus obtaining an existence result that rigorously recovers the Darcy model, which is formally derived by setting $\eta\equiv 0$ in system~\Ipbl,  as a mathematically consistent asymptotic regime.}
In Section~\ref{UNIQUENSS}, we analyze uniqueness and continuous dependence for the general system by assuming that the coefficients $\eta$ and $m$ are just positive constants.

%%%%%%%%%%%%%%%%%%%%%%%%%%%%%%%%%%%%%%%%%%%%%%%%%%%%%%%%%%%%%%%%%%%%%%%%%%%%%%
\section{Notation, assumptions and results}
\label{STATEMENT}
\setcounter{equation}{0}

\juerg{We assume that $\Omega\subset\erre^3$  is a bounded and connected open set}  with smooth boundary~$\Gamma$.
The symbol $|\Omega|$ denotes its Lebesgue measure, and $\dn$ stands for the outward normal derivative on~$\pier{\Gamma}$.
%Moreover, we keep the notation $Q:=\Omega\times(0,T)$.
For any Banach space~$X$, the symbols $\norma\cpto_X$ and $X^*$ indicate 
the corresponding norm and its dual space,
with \anold{a few} exceptions for the notation of the norm that are listed below.
\anold{The classical Lebesgue and Sobolev spaces on $\Omega$, corresponding to each $1 \leq p \leq \infty$ and $k \geq 0$, are denoted by $L^p(\Omega)$ and $W^{k,p}(\Omega)$, respectively. Their associated norms are written as $\norma{\cdot}_{L^p(\Omega)} = :\norma{\cdot}_p$ and $\norma{\cdot}_{W^{k,p}(\Omega)}$.}
Next, we~set
\begin{align}
  & H := \Ldue , \quad  
  V := \Huno , \quad
  W := \graffe{z\in\Hdue: \ \dn z=0 \ \hbox{on $\Gamma$}} , 
  \label{defspazia}
  \\
  & \HH := \juerg{H\times H\times H}, \quad
  \VV := \juerg{V\times V\times V}
  \aand
  \VVz := \graffe{\zz\in\VV : \ \div\zz=0}.
 \label{defspazib}
\end{align}
\Accorpa\Defspazi defspazia defspazib
Similarly, we use the boldface characters to denote powers of the Lebesgue and Sobolev spaces,
whose elements are vector valued functions.
For instance, we \anold{employ $\LLx p:= (\Lx p)^3$, for every $1 \leq p \leq \infty$, and}  $\HHx2=(\Hdue)^3$.
To~simplify  notation, the norms in the special cases $H$ and $\HH$ are \juerg{both} indicated by~$\norma\cpto$ 
\anold{,} without any subscript\anold{.}
The symbol $\norma\cpto_\infty$ might denote the norm in each of the spaces 
$\Linfty$, $\LQ\infty$ and $L^\infty(0,T)$, if no confusion can arise.
\anold{F}or simplicity, we use the same symbol for the norm in some space and the norm in any power of it,
as we have announced for $H$ and~$\HH$.
So, for example, the norms in $\VV$ and $\VVz$ are simply denoted by $\normaV\cpto$, 
since $\VV=\juerg{V\times V\times V}$ and $\VVz$ is a subspace of \anold{the latter}.

\anold{It \junew{is} worth recalling that} $V$ is dense in~$H$.
Therefore, we can adopt the usual framework of Hilbert triplets obtained by standard identifications.
Namely, we have that
\begin{align}
  & \< y,z > = \textstyle\iO y z 
  \quad \hbox{for every $y\in H$ and $z\in V$},
  \quad \hbox{so that} \quad
  V \emb H \emb \Vp .
  \label{embeddings}
\end{align}
In \eqref{embeddings}, the symbol $\< \cpto,\cpto >$ denotes the duality pairing between $\Vp$ and~$V$.
%The same notation is adopted throughout the paper.
Besides the space $\VVz$ already introduced, we also make use of the space
\Beq
  \HHz := \graffe{\zz\in\HH : \ \div\zz=0},
  \label{defHHz}
\Eeq
where the divergence is understood in the sense of distributions.
We notice at once that the embedding
\Beq
  \VVz \emb \HHz 
  \label{embeddingsbis}
\Eeq
is dense (see, e.g., \cite[Cor.~2.3]{Kato}).
We \anold{additionally} remark that all of the embeddings in \eqref{embeddings} and \eqref{embeddingsbis} are compact. 
Next, we recall the symbol $D\vv$ introduced in \eqref{symgrad} for the symmetrized gradient of the velocity~$\vv$,
whose use will be extended to any vector field $\zz\in\VV$.
Finally, we recall the standard notation \pier{for} the scalar product and the norm of matrices, namel\anold{y}
\Beq
  A:B := \somma {i,j}13 a_{ij} b_{ij}
  \aand
  |A|^2 := A:A 
  \quad \hbox{for $A=(a_{ij}),\,B=(b_{ij})\in\erre^{3\times3}$}\anold{.}
  \label{matrices}
\Eeq

\anold{We now list the assumptions regarding the structure of the system under consideration.}
As for the functions $\eta$, $\lambda$ and~$m$\anold{,} we assume that
\begin{align}
  & \hbox{$\eta:\erre\to\erre$ \ is \Lip\ continuous,}
  \non
  \\
  & \quad \hbox{and \ $\lambda,\,m:\erre\to\erre$ are locally \Lip\ continuous\pier{; all are positive, and}}
  \label{hpcoeff}
  \\
  & \anold{0<{}}\etamin \leq \eta(s) \leq \etamax \,, \quad
 \anold{0<{}}\lambdamin \leq \lambda(s) \leq \lambdamax
  \aand
  \anold{0<{}}\mmin \leq m(s) \leq \mmax,
  \non
  \\
  & \quad \hbox{for every $s\in\erre$ and some positive constants 
    \pier{$\etamin$, $\etamax$, $\lambdamin$, $\lambdamax$, $\mmin$ and $\mmax$}}.
  \label{hppositive}
\end{align}
Moreover, we postulate that
\an{
\begin{align}
  \non 
  & \nu \in \erre \, , \quad \pier{\sigma \in\erre,}  
  \aand
  S : \erre \to \erre
  \quad \hbox{is given by $S(\phi) = - \sigma \phi + h(\phi)$\,,}
  \\
  & \quad \hbox{\junew{where $h:\erre \to \erre$ is} bounded and \Lip\ continuous}.
  \label{hpSnu}
\end{align}
}
Finally, we assume that
\begin{align}
  & F : \erre \to \erre
  \quad \hbox{is of class \pier{$C^4$}, and \pier{$f:=F'$ denotes its derivative;}}
  \label{hppot}
  \\
  & \lim_{|s|\to+\infty} \frac{f(s)}s = + \infty\,\pier{;}
  \label{hpf}
  \\
  & f'(s) \geq -C_1 \,, \quad
  |F(s)| \leq C_2 \, (|s f(s)| + 1),
  \aand
  \pier{|sf'(s)| \leq C_3 (|f(s)| + 1),}
  \non
  \\
  & \quad \hbox{for every $s\in\erre$ and some positive constants $C_1$, $C_2$ and $C_3$}\,.
  \label{hpFf}
\end{align}
\Accorpa\Hpstruttura hpcoeff hpFf
Notice that the last inequality and the continuity of $f'$ imply that\pier{%
\Beq
|f' (s) |\leq C'_3(|f(s)|+1) \quad \hbox{for every $s \in \erre$,} 
\label{pier10}
\Eeq
for} a suitable constant~$C'_3$.

\anold{T}he conditions~\eqref{hppot}--\eqref{pier10} imposed on the potential \pier{can be compared with those in~\cite[assumptions~(2.4)--(2.9)]{CGSS-SIAM}, and are fulfilled, 
in particular, by the classical regular potential introduced in~\eqref{Freg}, as well as by any lambda-convex potential $F$ exhibiting polynomial growth.}

For the data\anold{,} we assume that
\begin{align}
  & \uu \in \L2\HH ,
  \label{regu}
  \\
  & \phiz \in W \,.
  \label{regphiz}
\end{align}
\Accorpa\Regdati regu regphiz

\last{Our main analytical objective is to establish the well-posedness of the system 
\Ipbl. To this end, we first highlight two key properties that allow equivalent formulations of some equations. As outlined in the Introduction, we then recast the velocity equation in a variational framework where the pressure variable $p$ no longer appears explicitly.
}

\last{To begin with, in place of equation \eqref{Iprima}, we adopt the variational formulation
\begin{align*}
	\iO \eta(\phi) D\vv : \nabla\zz
  + \iO \lambda(\phi) \vv \cdot \zz
  = \iO \mu \nabla\phi \cdot \zz
  + \iO \uu \cdot \zz
%  \,,
   \quad \hbox{for every $\zz\in\VVz$}\,,
\end{align*}
where the dot and colon denote the Euclidean inner product and the matrix inner product, respectively.
This reformulation is justified by the following observation concerning an equation of type \eqref{Iprima} 
for a fixed time and given $\badphi\in W$ and $\ff\in\HH$. Namely, from~\cite[Prop.~2.43]{E_thesis}} (see also \cite[Lemma~1.5]{EG1} and similar results in \cite{EL, EG_asy}), there exists a unique pair $(\yy,\pi)$ satisfying
\begin{align}
  & \yy \in \VVz
  \aand
  \pi \in H\,,
  \label{abte1}
  \\[2mm]
  & - \div\T(\badphi,\yy,\pi) + \lambda(\badphi) \yy = \ff
  \quad \hbox{in $\Omega$}\,,
  \label{abte2}
  \\[2mm]
  & \T(\badphi,\yy,\pi) \nn = \0 
  \quad \hbox{on $\Gamma$}\,.
  \label{abte3}
\end{align}
On the other hand, on account of \eqref{hppositive}, \juerg{it follows from the Lax--Milgram theorem} that
there exists a unique $\yy\in\VVz$ such that
\Beq
  \anold{a(\yy,\zz):=}\iO \bigl( \eta(\badphi) D\yy : \nabla\zz + \lambda(\badphi) \yy \cdot \zz \bigr) 
  = \iO \ff \cdot \zz
  \quad \hbox{for every $\zz\in\VVz$}\,. 
  \label{laxmil}
\Eeq
\juerg{Indeed, we} have that
\Beq
  \iO \bigl( \eta(\badphi) D\yy : \nabla\zz + \lambda(\badphi) \yy \cdot \zz \bigr) 
  \leq \max\graffe{\etamax,\lambdamax} \, \normaV\yy \, \normaV\zz
  \quad \hbox{for every $\yy,\zz\in\VVz$}\,, 
  \non
\Eeq
that is, the bilinear form \anold{$a:\VVz\times \VVz \to \erre$} given by the \lhs\ of \eqref{laxmil} is continuous on~$\VVz\times\VVz$.
Moreover, there hold the identity 
\Beq
  D\zz : \nabla\zz = |D\zz|^2
  \quad \hbox{for every $\zz\in\VV$}\anold{\,,}
  \label{identity}
\Eeq
and the Korn inequality
\Beq
  \normaV\zz^2 \leq \CK \iO (|D\zz|^2 + |\zz|^2)
  \quad \hbox{for every $\zz\in\VV$} \,,
  \label{korn}
\Eeq
with some constant $\CK$ depending only on~$\Omega$\juerg{. Hence,} \pier{by setting 
$\alpha := \min \graffe{\etamin,\lambdamin}/\CK$}, \anold{it follows that}
\begin{align}
  & \anold{a(\zz,\zz)=}\iO \bigl( \eta(\badphi) \pier{| D\zz|^2} + \lambda(\badphi) |\zz|^2 \bigr) 
  \geq \alpha \, \normaV\zz^2
  \quad \hbox{for every $\zz\in\VVz$}\,.
  \label{coercive}
%  \\
%  & \hbox{where} \quad 
%  \alpha := \frac{\min\graffe{\etamin,\lambdamin}}\CK\,.
%  \label{defalpha}
\end{align}
\juerg{Therefore}, the bilinear form under consideration is also coercive.
Finally, the component $\yy$ of the solution $(\yy,\pi)$ to \accorpa{abte1}{abte3} 
belongs to $\VVz$ and also satisfies \eqref{laxmil}, thanks to the identities
\begin{align}
  & \iO \eta(\badphi) D\yy : \nabla\zz
  = \iO \bigl( \eta(\badphi) D\yy - \pi \I \bigr) : \nabla\zz
  \non
  \\
  & = \iO \anold{\T}(\badphi,\yy,\pi) : \nabla\zz
  = - \iO \div\T(\badphi,\yy,\pi) \cdot \zz \,,
  \non
\end{align}
which\juerg{, owing to \eqref{abte3}, hold true for every $\zz\in\VVz$}.
\anold{Thus,} it is equivalent to look either for the unique solution $\yy\in\VVz$ to \eqref{laxmil}
or for the first component of the solution $(\yy,\pi)$ to \accorpa{abte1}{abte3}.
Since the component $\phi$ of the solution to the problem we want to consider is required to belong to $\L\infty W$
and the \rhs\ of \eqref{Iprima} is expected to belong to $\L2\HH$,
we can apply the above \anold{correspondence} to $\badphi=\phi(t)$ \aat\
and replace \eqref{Iprima} by the forthcoming variational equation~\eqref{prima}.

\anold{We now consider an equivalent formulation of the problem.
When appropriate, we eliminate $w$ by substituting the expression provided by \eqref{Iquarta}
 \junew{in} \eqref{Iterza}. This leads to the following equation:}
\Beq
  - \Delta \bigl( -\Delta\phi + f(\phi) \bigr)
  + \bigl( f'(\phi) + \nu \bigr) \bigl( -\Delta\phi + f(\phi) \bigr)
  = \mu\anold{,}
  \label{quintastrong}
\Eeq
and can recover $w$ and \eqref{Iterza} by using \eqref{Iquarta} as a definition of~$w$.
However, we will write all the equations in their variational forms, for convenience.

Here is the precise formulation of the problem under investigation.
We look for a quadruplet $\soluz$ with the regularity
\begin{align}
  & \vv \in \L2\VVz\,,
  \label{regv}
  \\
  & \phi \in \H1\Vp \cap \L\infty W \cap \anold{\L2{\Hx5}}\,,
  \label{regphi}
  \\
  & \mu \in \L2V\,,
  \label{regmu}
  \\
  & w \in \L\infty H \cap \anold{\L2 {\Hx3 \cap W}}\,,
  \label{regw}
\end{align}
\Accorpa\Regsoluz regv regw
that solves the variational equations
\begin{align}
  & \iO \bigl(
    \eta(\phi) D\vv : \nabla\zz
    + \lambda(\phi) \vv \cdot \zz
  \bigr)
  = \iO (\mu \nabla\phi + \uu) \cdot \zz
  \non
  \\
  & \quad \hbox{for every $\zz\in\VVz$ and \aet}\,,
  \label{prima}
  \\
  & \< \dt\phi , z > 
  + \iO \vv \cdot \nabla\phi \, z
  + \iO m(\phi) \nabla\mu \cdot \nabla z
%  \pier{{}+ \sigma \iO \phi z}
  = \iO S(\phi) z
  \non
  \\
  & \quad \hbox{for every $z\in V$ and \aet}\,,
  \label{seconda}
  \\
  & 
  \anold{\iO \nabla w 	\cdot \nabla  z}
%  \iO w (-\Delta z)
  \juerg{\,+} \iO \bigl( f'(\phi) + \nu \bigr) w z 
  = \iO \mu z
  \non
  \\
  & \quad \hbox{for every $z\in \anold{V}$ and \aet}\,,
  \label{terza}
  \\
  & \iO \nabla\phi \cdot \nabla z
  + \iO f(\phi) z
  = \iO w z
  \non
  \\
  & \quad \hbox{for every $z\in V$ and \aet}\,,
  \label{quarta}
\end{align}
as well as the initial condition
\Beq
  \phi(0) = \phiz \,.
  \label{cauchy}
\Eeq
\Accorpa\Pbl prima cauchy
Notice that \eqref{terza} and \eqref{quarta}
\anold{contain the} homogeneous Neumann boundary condition for $w$ and $\phi$ in a generalized sense.
However, these equations can be replaced by their strong forms \eqref{Iterza} and \eqref{Iquarta},
thanks to the regularity of $w$ and $\phi$ required in \eqref{regw} and~\eqref{regphi}\anold{,} which also contain the homogeneous Neumann boundary condition\pier{s.}
By taking this into account, we can equivalently replace \pier{\accorpa{terza}{quarta}}~by
\begin{align}
  & 
  \anold{\iO \nabla \bigl( 
    - \Delta\phi + f(\phi)
  \bigr) \cdot \nabla z}
%  \iO \bigl( 
%    - \Delta\phi + f(\phi)
%  \bigr) (-\Delta z)
  + \iO \bigl( f'(\phi) + \nu \bigr) \bigl( - \Delta\phi + f(\phi) \bigr) z
  = \iO \mu z
  \non
  \\
  & \quad \hbox{for every $z\in \anold{V}$ and \aet},
  \label{quinta}
\end{align}
which is the weak form of \eqref{quintastrong}, and keep \eqref{quarta}\anold{,} or \eqref{Iquarta}\anold{,} as a definition of~$w$.

Here is our first result.

\Bthm
\label{Wellposedness}
\pier{Under the assumptions~\Hpstruttura\ on the structure and~\Regdati\ on the data,
there} exists at least one quadruplet $\soluz$ \pier{that fulfills the regularity requirements \Regsoluz, solves} Problem \Pbl,  and satisfies the estimate
\begin{align}
  & \norma\vv_{\L2\VVz}
  + \norma\phi_{\H1\Vp\cap\L\infty W\cap\anold{\L2{\Hx5}}}
  \non
  \\
  & \quad {}
  + \norma\mu_{\L2V}
  + \norma w_{\L\infty H\cap\anold{\L2{\Hx3 \cap W}}}
  \leq K_1\,,
  \label{stability}
\end{align}
with a constant $K_1$ that depends only on the structure of the system, $\Omega$, $T$ 
and an upper bound for the norms of the data related to \Regdati.
Moreover, the solution is unique if both $\eta$ and $m$ are positive constants.
\Ethm

\anold{%
\pier{We now establish an existence result for a system closely related to~\Pbl, 
but governed by the so-called Darcy law,} which formally corresponds to 
setting the shear viscosity $\eta \equiv 0$. In this asymptotic framework, 
particular care must be devoted to the boundary condition  $\T(\phi, \vv, p)\nn = \mathbf{0}$, which, as is well known \pier{(see, e.g., \cite{EG_asy} and also \cite{KS2})}, \an{corresponds} to the boundary condition $p = 0$ on $\Sigma$ \an{in the vanishing viscosity limit}.
\begin{theorem}\label{Darcy}
Suppose that \pier{the assumptions~\Hpstruttura\ on the structure and~\Regdati\ on the data are fulfilled.}
Furthermore, let $\{\eta_n\}_{n \in \mathbb{N}}$  be  a sequence of viscosity functions such that, for each fixed $n \in \mathbb{N}$, $\eta_n$ is compatible with \junew{the} assumptions \eqref{hpcoeff} and \eqref{hppositive}\pier{, for two sequences of positive values $ {\etamin}_{,n} $ and ${\etamax}_{,n} $ 
\junew{satisfying} $\ 0<  {\etamin}_{,n} \leq \eta_n(s) \leq {\etamax}_{,n} \ $ for every $s\in \erre$.}  
We further assume that
\begin{align*}
\pier{\|\eta_n\|_{\pier{L^\infty(\mathbb{R})}} \, \hbox{ tends to } \, 0 \,  \text{ as } n \to \infty.}
\end{align*}
For any $n \in \mathbb{N}$, let $(\vvn, \varphi_n, \mu_n, w_n)$ denote \pier{a solution \junew{to} Problem~\Pbl\ as found in} Theorem~\ref{Wellposedness} associated with the \an{function $\eta_n$}.
% and initial datum $\phiz$.
\junew{Then there is a subsequence of the sequence $\{(\vvn,  \varphi_n, \mu_n,w_n)\}$, which is again labeled by
$n\in\enne$, that converges to a solution $(\vv,  \varphi, \mu,w)$ of the corresponding Darcy system in the following sense:}  as $n\to \infty$, \pier{it holds that} 
\begin{align}
	\vvn  \to \vv &\quad \text{weakly in $\L2 {\HH_0}$}\,,   \label{convD1}
	\\
	\eta_n(\phin) D \vvn
	 \to \0 &\quad \text{\pier{strongly} in $L^2(0,T;\HH^3)$}\,,  \label{convD2}
	\\[1mm]
%		p_n  \to p &\quad \text{weakly in $\L2 H$}\,,
%	\\
		\phi_n  \to \phi &\quad \text{weakly star in }\non\\
		&\quad\quad \H1\Vp \cap \L\infty W \cap  \L2 {\Hx5}\,, \label{convD3}
	\\[1mm]
		\mu_n \to \mu &\quad \text{weakly in $\L2 V$}\,, \label{convD4}
		\\
		w_n \to w &\quad \text{weakly star in $\L\infty H \cap \L2 {\Hx3 \cap W}$}\,, \label{convD5}
\end{align}
\pier{where $(\vv,  \varphi, \mu,w)$ solves} the variational \pier{equations}
\begin{align}
	& \iO 
    \lambda(\phi) \vv \cdot \zz
  = \iO (\mu \nabla\phi + \uu) \cdot \zz \non\\
   & \quad  \hbox{for every \pier{$\zz\in\HHz$} and \aet}\,,
   \label{primaD}
   \\[2mm]
   & \< \dt\phi , z > \pier{{}+ \< \vv \cdot \nabla \phi ,z> }
%  - \iO (\vv \phi) \cdot \nabla z
  + \iO m(\phi) \nabla\mu \cdot \nabla z
%    \pier{{}+ \sigma \iO \phi z}
  = \iO S(\phi) z   \non\\
  &\quad \hbox{for every $z\in V$ and \aet}\,,
  \label{secondaD}
\end{align}
\pier{as well as~\eqref{terza}--\eqref{quarta} and the} initial condition \eqref{quinta}.
% and the boundary condition $p=0$ \aeS.
\end{theorem}}
\pier{We point out that the equations \eqref{primaD} and \eqref{secondaD} are \an{consistently formulated}. 
Indeed, owing to \eqref{regphi} and \eqref{regmu}, the integral on the right-hand side 
of \eqref{primaD} involves the term 
\[
  \mu \nabla \phi,
\]
which is the product of $\mu \in L^2\big(0,T;L^4(\Omega)\big)$ and 
$\nabla \phi \in L^\infty\big(0,T;\LLx4\big)$. 
Hence, one has $\mu \nabla \phi \in L^2(0,T;\HH)$. 
On the other hand, the second term in \eqref{secondaD} contains 
\[
  \vv \cdot \nabla \phi,
\]
where $\vv $ is just in $ L^2(0,T;\HH)$\junew{. But, in fact, we have} 
$\vv \cdot \nabla \phi \in L^2\big(0,T;L^{4/3}(\Omega)\big)$. 
Since $L^{4/3}(\Omega)$ is continuously embedded \junew{in} $\Vp$, 
we deduce that $\vv \cdot \nabla \phi \in L^2(0,T;\Vp)$.}
\an{Let us also remark that, due to the above regularities, it even holds  that $\vv \cdot \nabla \phi \in \L1 H$ so that the corresponding term in \eqref{secondaD} may be also written as an integral.
}

\pier{We now return to Problem~\Pbl\ and address the issues of uniqueness 
and continuous dependence of the solution with respect to the forcing term 
$\uu$ appearing in the first equation.}

\Bthm
\label{Contdep}
Assume that \Hpstruttura\ \juerg{are fulfilled and} that $\eta$ and $m$ are positive constants. Moreover,  assume \eqref{regphiz} for the initial datum. \juerg{Whenever}  $\uu_i$, $i=1,2$, belong to $\L2\HH$ and $(\vv_i,\phi_i,\mu_i,w_i)$ are the corresponding solutions,
then the estimate
\begin{align}
  & \norma\vv_{\L2\VV}
  + \norma\phi_{\C0V\cap\L2{\Hx4}}
  + \norma\mu_{\L2H}
  \non
  \\
  & \quad {}
  + \norma w_{\L2W}
  \leq K_2 \, \norma\uu_{\L2\HH}
  \label{contdep}
\end{align}
\juerg{holds true} for the differences $\soluz=(\vv_1,\phi_1,\mu_1,w_1)-(\vv_2,\phi_2,\mu_2,w_2)$ and $\uu=\uu_1-\uu_2$,
with a constant $K_2$ that depends only on the structure of the system, $\Omega$, $T$, the initial datum~$\phiz$
and an upper bound for the norms of $\uu_1$ and $\uu_2$ in~$\L2\HH$.
\Ethm

The remainder of the section is devoted to the collection of some useful tools.
First of all, besides the \Holder\ and Cauchy--Schwarz inequalities
and the Korn inequality \eqref{korn}, we \juerg{widely use Young's} inequality
\Beq
  ab \leq \delta a^2 + \frac 1{4\delta} \, b^2
  \quad \hbox{for every $a,b\in\erre$ and $\delta>0$}.
  \label{young}
\Eeq
We also account for the Sobolev and Poincar\'e inequalities, 
as well as for some inequalities associated \juerg{with} the elliptic regularity theory 
and the compact embeddings between Sobolev spaces (via Ehrling's lemma).
\anold{More precisely, the following estimates hold:}
\begin{align}
  & \norma z_q
  \leq \CS \, \normaV z
  \quad \hbox{for every $z\in V$ and $q\in[1,6]$}.
  \label{sobolev1}
  \\[2mm]
  \separa
  & \norma z_\infty
  \leq \CS \, \normaW z
  \quad \hbox{for every $z\in \anold{W}$}.
  \label{sobolev2}
  \\[2mm]
  \separa
  & \normaV z
  \leq \CP \, \bigl( \norma{\nabla z} + |\zbar| \bigr)
  \quad \hbox{for every $z\in V$}.
  \label{poincare}
  \\[2mm]
  \separa
  & \normaW z
  \leq \CE \, \bigl( \norma{\Delta z} + \norma z \bigr)
  \label{elliptic1}
  \quad \hbox{for every $z\in W$}.
  \\[2mm]
  & \norma z_{\Hx4}
  \leq \CE \, \bigl( \norma{\Delta^2 z} + \norma z \bigr)
  \quad \hbox{for every $z\in\Hx4$ with $z,\Delta z\in W$}.
  \label{elliptic2}
  \\[2mm]
  \separa
  & \normaV z 
  \leq \delta \, \norma{\Delta v} + \Cdelta \, \norma z  
  \quad \hbox{for every $z\in W$ \anold{and every $\delta>0$}}.
  \label{compact1}
  \\[2mm]
  & \norma z_{\Hx3}
  \leq \delta \, \norma{\Delta^2 z} + \Cdelta \, \normaV z  
  \non
  \\
  & \quad \hbox{for every $z\in\Hx4$ with $z,\Delta z\in W$ and every $\delta>0$}.
  \label{compact2}
\end{align}
Here, $\zbar$ in \eqref{poincare} denotes the mean value of~$z$. 
In general, we~set
\Beq
  \zbar := \frac 1{|\Omega|} \, \iO z
  \quad \hbox{for $z\in\Luno$},
  \label{mean}
\Eeq
and we use the same notation for time-dependent functions.
In the above inequalities \accorpa{sobolev1}{elliptic2}, 
the constants on the \rhs s  depend only on~$\Omega$, 
while $\Cdelta$ in \accorpa{compact1}{compact2} depends on both $\Omega$ and~$\delta$.
\anold{Analogous inequalities naturally extend to vector-valued functions as well without the need for repetition.}
%For instance, we can replace $\Lx q$ and $V$ in \eqref{sobolev1} 
%by $\LLx q$ and~$\VV$, respectively.
In connection with \eqref{elliptic1} and \eqref{elliptic2}, 
we recall that $z$ belongs to $W$ whenever $z\in V$, $\Delta z\in H$ 
and the homogeneous Neumann boundary condition is satisfied in the usual weak sense,
and that $z$ belongs to $\Hx4$ whenever both $z$ and $\Delta z$ belong to~$W$.

We conclude this section by introducing a convention that simplifies \anold{some of} the \juerg{subsequent calculations}.
The lowercase symbol $c$ denotes a generic constant
that depends only on $\Omega$, $T$, the structure of the system, and an upper bound for the norms of the data given in~\Regdati. 
In particular, $c$~is independent of the parameter $n$ that we introduce in the next section.
Notice that the value of $c$ may change from line to line and even within the same line \anold{of computations}.
Furthermore, whenever some positive constants, e.g., $\delta$ or~$M$, occur in a computation,
we adopt the corresponding subscript by writing, e.g., $\cdelta$ and~$\cM$, instead of a general~$c$, 
\anold{to emphasize} that these constants also depend on the parameter under consideration. 
On the contrary, specific constants referenced explicitly are denoted by different symbols,
like in \eqref{stability} and~\eqref{sobolev1}, where different characters are used.

%%%%%%%%%%%%%%%%%%%%%%%%%%%%%%%%%%%%%%%%%%%%%%%%%%%%%%%%%%%%%%%%%%%%%%%%%%%%%%

\section{Existence \pier{of solutions}}
\label{EXISTENCE}
\setcounter{equation}{0}

This section is devoted to prove the existence part of Theorem~\ref{Wellposedness}.
More precisely, we prove the existence of a quadruplet $\soluz$ 
with the regularity specified in \Regsoluz\ that solves Problem \Pbl\
and satisfies the inequality \eqref{stability} with a constant $K_1$ \juerg{having the properties stated in the assertion}.
Our \anold{rigorous} argument is based on a Faedo--Galerkin scheme 
associated with the introduction of a viscosity term in the first equation \eqref{prima}.

To this end, we introduce the sequence $\graffe\lambdaj_{j\geq1}$ of the eigenvalues 
and an orthonormal system $\graffe\ej_{j\geq1}$ of corresponding eigenfunctions
of the Neumann problem for the Laplace equation, \anold{that is,}
\begin{align}
  & 0 = \lambda_1 < \lambda_2 \leq \lambda_3 \leq \dots
  \aand \lim_{j\to\infty} \lambdaj = + \infty\,,
  \label{eigenvalues}
  \\
  & \ej\in V
  \aand
  \iO \nabla\ej \cdot \nabla \anold{z}
  = \lambda_j \iO \ej \anold{z}
  \quad \hbox{for every $\anold{z}\in V$ and $j=1,2,\dots,$}
  \label{eigenfunctions}
  \\
  & \iO \ei \ej = \delta_{ij}
  \quad \hbox{for $i,j=1,2,\dots,$}
  \aand 
  \hbox{$\graffe\ej_{j\geq1}$ is a complete system in $H$,}
  \label{orthogonality}
\end{align}
where $\delta_{ij}$ is the Kronecker symbol.
Similarly, we consider the eigenvalue problem
\Beq
  \eej \in \VVz
  \aand
  \iO D\eej : \nabla\zz = \lambdazj \iO \eej \cdot \zz\,,
  \quad \hbox{for every $\zz\in\VVz$}.
  \label{eeigenfunctions}
\Eeq
Since the bilinear form appearing on the \lhs\ of \eqref{eeigenfunctions} is continuous and weakly coercive on $\VVz\times\VVz$
(see \eqref{coercive} with $\eta$ and $\lambda$ replaced by~$1$),
and since $\VVz$ is compactly embedded in~$\HHz$,
the general \anold{abstract} theory can be applied, and the above problem provides
the sequence $\graffe\lambdazj_{j\geq1}$ of the eigenvalues 
and an orthonormal system $\graffe\eej_{j\geq1}$ of corresponding eigenfunctions.
We thus have
\begin{align}
  & 0 = \lambda_{0,1} \leq \lambda_{0,2} \leq \lambda_{0,3} \leq \dots
  \aand \lim_{j\to\infty} \lambdazj = + \infty\,,
  \label{eeigenvalues}
  \\
  & \iO \eei \cdot \eej = \delta_{ij}
  \quad \hbox{for $i,j=1,2,\dots$},
  \aand 
  \hbox{$\graffe\eej_{j\geq1}$ is complete in $\HHz$} \,.
  \label{oorthogonality}
\end{align}
We \anold{then}~set
\Beq
  \Vn := \Span \graffe{e_1,\dots,e_n}
  \aand
  \VVzn := \Span \graffe{\ee_{0,1},\dots,\ee_{0,n}},
  \quad \hbox{for $n=1,2,\dots,$}
  \label{defVn}
\Eeq
and observe that each $\Vn$ is included in \juerg{$W$} and that $\Delta z$ belongs to $\Vn$ for every $z\in\Vn$.
Moreover, the constant functions belong to~$V_1$ and thus to every~$\Vn$.
Finally, we stress that the union of the spaces $\Vn$ is dense in both $V$ and~$H$.
Similarly, the union of the spaces $\VVzn$ is dense in both $\VVz$ and~$\HHz$.

\Brem
\label{Projection}
Let us make some observations regarding the orthogonal projection operator $\Pn:H\to\Vn$\anold{, where  the orthogonality is understood with respect to the standard inner product of~$H$.}
Let $Y$ \anold{denote} any of the spaces $H$, $V$, \anold{or}~$W$.
As shown in \cite[Rem~3.3]{CGSS-SIAM} and in \cite[Rem.~3.1]{CG2}, we have\anold{,} for every $z\in Y$\anold{,} that 
\Beq
  \norma{\Pn z}_Y \leq \CO \, \norma z_Y
  \aand
  \Pn z \to z
  \quad \hbox{strongly in $Y$},
  \label{convPn}
\Eeq
with a constant $\CO$ that depends only on~$\Omega$.
%Moreover, one can take $\CO=1$ provided that the norm in $Y$ is well chosen  \juerg{\bf I do not understand the rest of this sentence, since
%the strong convergence was just stated}
%so that the above convergence actually is strong by uniform convexity.
Next, if the operator $\Pn$ is extended to spaces of time-dependent functions
(i.e.,~for $z\in\L2Y$ the function $\zn$ is defined by
$\zn(t):=\Pn(z(t))$ \aat),
then
\Beq
  \norma\zn_{\L2Y} \leq \pier{\CO} \norma z_{\L2Y},
  \aand
  \zn \to z 
  \quad \hbox{strongly in $\L2Y$}.
  \label{convPnbis}
\Eeq
These observations are useful when we let $n$ tend to infinity in the Faedo\anold{--}Galerkin scheme introduced below.
Unfortunately, it is not clear whether similar inequalities and strong convergence properties
hold for the analogous projection operator from $\HHz$ onto~$\VVzn$.
\Erem

At this point, we are ready to introduce the Faedo--Galerkin scheme mentioned above.
It consists in looking for a triplet
\Beq
  \soluzn \in \H1\VVzn \times \H1\Vn \times \L2\Vn
  \label{regsoluzn}
\Eeq
that solves the variational equations 
\begin{align}
  & \frac 1n \iO \dt\vvn \cdot \zz 
  + \iO \bigl(
    \eta(\phin) D\vvn : \nabla\zz
    + \lambda(\phin) \vvn \cdot \zz
  \bigr)
  = \iO (\mun \nabla\phin + \uu) \cdot \zz
  \non
  \\
  & \quad \hbox{for every $\zz\in\VVzn$ and \aet},
  \label{priman}
  \\
  & \iO \dt\phin \, z 
  + \iO \vvn \cdot \nabla\phin \, z
  + \iO m(\phin) \nabla\mun \cdot \nabla z
%  \pier{{}+ \sigma \iO \phin z} 
= \iO S(\phin) z
  \non
  \\
  & \quad \hbox{for every $z\in\Vn$ and \aet},
  \label{secondan}
  \\
  & \iO \bigl( 
    - \Delta\phin + f(\phin)
  \bigr) (-\Delta z)
  + \iO \bigl( f'(\phin) + \nu \bigr) \bigl( - \Delta\phin + f(\phin) \bigr) z
  = \iO \mun z
  \non
  \\
  & \quad \hbox{for every $z\in\Vn$ and \aet} ,
  \label{quintan}
\end{align}
and satisfies the initial conditions 
\Beq
  \vvn(0) = \0
  \aand
  \phi(0) = \Pn\phiz \,.
  \label{cauchyn}
\Eeq
\Accorpa\Pbln priman cauchyn

First of all, we have to prove that the above problem is well posed.

\step
Well-posedness of the discrete problem

\anold{We show that the problem admits a unique maximal solution. Then, by proving several a priori estimates, we demonstrate that this maximal solution extends globally. This aspect will be revisited later.}

\anold{We express the unknowns as expansions} in terms of the eigenfunctions $\eej$ and $\ej$ as follows:
\begin{align}
  & \vvn(t) = \somma j1n \vnj(t) \eej \,, \quad
  \phin(t) = \somma j1n \phinj(t) \ej \,,
  \non
  \\
  & \aand
  \mun(t) = \somma j1n \munj(t) \ej\,,
  \quad \aat,
  \non
\end{align}
where the coefficients are required to satisfy
\Beq
  \vnj \in H^1(0,T) , \quad
  \phinj \in H^1(0,T) 
  \aand
  \munj \in L^2(0,T) \,.
  \non
\Eeq
Then, we insert the above expansions in the variational equations \anold{\Pbln}.
Since it suffices to choose $\zz=\eei$ in \eqref{priman} and $\zeta=\ei$ in \eqref{secondan} and \eqref{quintan}, for $i=1,2,\dots,n$,
\juerg{we obtain from the orthogonality properties of the eigenfunctions}  a system of the form
\begin{align}
  & \frac 1n \, \vni'
  = \calF_1 \bigl( (\vnj)_{j=1}^n , (\phinj)_{j=1}^n ,(\munj)_{j=1}^n \bigr),
  \label{primanbis}
  \\
  & \phini'
  = \calF_2 \bigl( (\vnj)_{j=1}^n , (\phinj)_{j=1}^n ,(\munj)_{j=1}^n \bigr),
  \label{secondanbis}
  \\
  & \muni
  =\pier{\calF_3 \bigl( (\phinj)_{j=1}^n \bigr)},   %\bigl( (\vnj)_{j=1}^n , (\phinj)_{j=1}^n ,(\munj)_{j=1}^n \bigr)
  \label{quintanbis}
\end{align}
all for $i=1,\dots,n$ and a.e.\ with respect to time,
\anold{the functions $\calF_1$, $\calF_2$ and $\calF_3$, naturally defined,} are locally \Lip\ continuous functions on $\erren\times\erren\times\erren$, 
thanks to \eqref{hpcoeff}, \eqref{hpSnu} and~\eqref{hppot}.
\anold{By multiplying}~\eqref{primanbis} by $n$ and \pier{using}~\eqref{quintanbis} to eliminate every $\munj$ in
\juerg{\eqref{primanbis} and} \eqref{secondanbis},
\pier{the problem} is reduced to a system of explicit ordinary differential equations in standard form,
ruled by locally \Lip\ continuous functions,
in the unknowns $(\vnj)_{j=1}^n$ and~$(\phinj)_{j=1}^n$.
On the other hand, \eqref{cauchyn} provide proper initial values.
Therefore, the Cauchy problem we obtain has a unique maximal solution ($(\vnj)_{j=1}^n,(\phinj)_{j=1}^n)$ \anold{defined in the interval~$[0,T_n)$ for some $T_n\in(0,T]$.}
Since $(\munj)_{\juerg{j=1}}^n$ can be recovered from~\eqref{quintanbis},
we conclude that Problem \Pbln\ has a unique maximal solution defined in the interval~$[0,T_n)$\anold{.}

\medskip

Now, we start estimating.
As said before, \anold{the estimates we establish also guarantee that} the solution $\soluzn$ to Problem \Pbln\ is global,
i.e., that $T_n=T$.
For this reason and in order to simplify the notation, we write $T$ instead of $T_n$ from now on.
\anold{Furthermore}, we avoid writing the subscript $n$ in performing calculations
and restore the \anold{rigorous} notation $\soluzn$ only at the end of each step.
\anold{Additionally}, for any given time-dependent test function $z$,
it is understood that the \juerg{respective} equations are always written at the time $t$ and then tested by $z(t)$ \aat, 
even though we avoid writing the time, for simplicity.

\step
First a priori estimate

We test \pier{equation}~\eqref{priman} \anold{with~}$\vvn$, \pier{equation}~\eqref{secondan} \anold{with} $\mun $, and \pier{equation}~\eqref{quintan} \anold{with~}$\dt\phin \pier{{}+\sigma \phin}$.
On account of~\eqref{identity}, we obtain~that 
%(where the subscript $n$ is not written)
\begin{align}
  & \frac 1{2n} \, \ddt \, \iO |\vv|^2
  + \iO \bigl( \eta(\phi) |D\vv|^2 + \lambda(\phi) |\vv|^2 \bigr)
  = \iO \mu \nabla\phi \cdot \vv
  + \iO \uu \cdot \vv\,,
  \non
  \\[2mm]
  & \iO \dt\phi \, \mu 
  + \iO \vv \cdot \nabla\phi \, \mu
  + \iO m(\phi) |\nabla\mu|^2
%   \pier{{} + \sigma\iO \phi \, \mu} 
  = \iO S(\phi) \mu\,,
  \non
  \\[2mm]
  & 
%	\anold{\iO \nabla \bigl( 
%    - \Delta\phi + f(\phi)
%  \bigr) \cdot \nabla  (\dt\phi)
%}	
\pier{	\iO \bigl( - \Delta\phi + f(\phi)
\bigr) (-\Delta\dt\phi)}
  + \iO \bigl( f'(\phi) + \nu \bigr) \bigl( - \Delta\phi + f(\phi) \bigr) \dt\phi
   \non
  \\
  &{} = 
  \pier{	- \sigma \iO \bigl( - \Delta\phi + f(\phi)
\bigr) (-\Delta \phi)
  -  \sigma \iO \bigl( f'(\phi) + \nu \bigr) \bigl( - \Delta\phi + f(\phi) \bigr) \phi}
    \non
  \\
  &\quad\pier{{}+ \iO \mu \, (\dt\phi  \pier{{}+\sigma \junew{\phi}}) } \,.
  \non
\end{align}
%\comm{[Pier: removed a modification by Andrea on the first term of the third line above since it was not appropriate (better to leave it as in~\eqref{quintan} and the equality below)]}\\
We add these identities to each other.
%\anold{, and owe to integration by parts}.
Since \pier{some} cancellations occur, we have~that
\begin{align}
  & \frac 1{2n} \, \ddt \, \iO |\vv|^2
  + \iO \bigl( \eta(\phi) |D\vv|^2 + \lambda(\phi) |\vv|^2 \bigr)
  + \iO m(\phi) |\nabla\mu|^2
  \non
  \\
  & \quad {}
  + \iO \bigl( 
    - \Delta\phi + f(\phi)
  \bigr) (-\Delta\dt\phi)
  + \iO \bigl( f'(\phi) + \nu \bigr) \bigl( - \Delta\phi + f(\phi) \bigr) \dt\phi
  \non
  \\
  & = \iO \uu \cdot \vv
  + \iO \an{h(\phi)} \mu 
  \pier{{}	- \sigma \iO \bigl( - \Delta\phi + f(\phi)
\bigr) \bigl(-\Delta \phi + f'(\phi) \phi + \nu \phi \bigr)}
  \,.
  \non
\end{align}
Now, we notice that the time derivative of the energy defined in \eqref{defE} is given~by
\begin{align}
  & \ddt \, \calE(\phi)
  = \iO \bigl(
    -\Delta\phi + f(\phi) 
  \bigr)
  \bigl(
    -\Delta\dt\phi + f'(\phi)\dt\phi 
  \bigr)  
  \non
  \\
  & \quad {}
  + \nu \iO \bigl( 
    \nabla\phi \cdot \nabla\dt\phi + f(\phi) \dt\phi 
  \bigr)
  \label{dtenergy}
\end{align}
and recall the coercivity inequality \eqref{coercive} \anold{along} with the definition \pier{of~$\alpha$ (see \eqref{coercive})} and the assumptions~\pier{\eqref{hppositive}.}
We deduce the inequality
\begin{align}
  &\frac 1{2n} \, \ddt \, \pier{\norma\vv^2}
  + \pier{\frac\alpha 2} \, \normaV\vv^2
  \pier{{}+ \frac12 \,\iO \eta(\phi) \pier{| D\vv|^2}}
  \pier{{}+ \frac{\lambda_*}2 \, \norma\vv^2}
  + \mmin \iO |\nabla\mu|^2
  + \ddt \, \calE(\phi)
  \nonumber
  \\
  &{}
  \leq \iO \uu \cdot \vv
  + \iO \an{h(\phi)}\mu \pier{{}	- \sigma \iO \bigl( - \Delta\phi + f(\phi)
\bigr) \bigl(-\Delta \phi + f'(\phi) \phi + \nu \phi \bigr)\, ,}
  \label{pier1}
\end{align}
\pier{written in some complicated form, which is however good for later purposes.}
At this point, we integrate over $(0,t)$ with respect to time\anold{, for any $t \in (0,T]$}.
However, since a part of the energy is multiplied by $\nu$ and this parameter is not required to be positive,
it is convenient to rewrite $\calE(\phi)$ in a different form by splitting the first integral into two parts
and making explicit calculations in the second one.
\pier{In the light of~\eqref{defE}, we} have
\begin{align}
  & \calE(\phi)
  = \frac14  \iO \bigl( -\Delta\phi + f(\phi) \bigr)^2
  + \frac14 \iO \bigl(
    |\Delta\phi|^2 + |f(\phi)|^2 
  \bigr)
  \non
  \\
  & \quad {}
  + \frac 12 \iO f'(\phi) |\nabla\phi|^2 
  + \nu \iO \Bigl( \frac 12 \, |\nabla\phi|^2 + F(\phi) \Bigr)\,.
  \label{energy}
\end{align}
Therefore, the \pier{time integration of~\eqref{pier1}} and some rearrangement yield, for every $t\in(0,T]$,
\begin{align}
  & \frac 1{2n} \, \pier{\norma{\vv(t)}^2}
  + \pier{\frac\alpha 2} \iot \normaV{\vv(s)}^2 \, ds
   \pier{{}+ \frac12 \,\intQt \eta(\phi) \pier{| D\vv|^2}}
  \pier{{}+ \frac{\lambda_*}2 \iot \norma{\vv(s)}^2 \, ds}
  \non
  \\
  & \quad {}
  + \mmin \intQt |\nabla\mu|^2
  + \frac14  \iO |\bigl( -\Delta\phi + f(\phi) \bigr)(t)|^2
  + \frac14 \iO \bigl(
    |\Delta\phi(t)|^2 + |f(\phi(t))|^2 
  \bigr)
  \non
  \\
  & {} \pier{{}\leq{}} - \frac 12 \iO f'(\phi(t)) |\nabla\phi(t)|^2 
  - \nu \iO \Bigl( \frac 12 \, |\nabla\phi(t)|^2 + F(\phi(t)) \Bigr)
  \non
  \\
  & \quad {} + \calE \pier{(\Pn \phiz)}
  + \intQt \uu \cdot \vv
  + \intQt \an{h}(\phi) \mu
  \non
  \\
  & \quad {}\pier{{}	- \sigma \intQt \bigl( - \Delta\phi + f(\phi)
\bigr) \bigl(-\Delta \phi + f'(\phi) \phi + \nu \phi \bigr)}
  \,,
  \label{perprimastima}
\end{align}
and we have to estimate the terms on the \rhs.
On account of the first assumption in \eqref{hpFf} and  of the compactness inequality \eqref{compact1},
we have~that
\begin{align}
  & - \frac 12 \iO f'(\phi(t)) |\nabla\phi(t)|^2 
  - \frac \nu 2 \iO |\nabla\phi(t)|^2
  \leq \frac {C_1 + |\nu|}2 \iO |\nabla\phi(t)|^2
  \non
  \\
  & \leq \frac 18 \iO |\Delta\phi(t)|^2
  + C' \iO |\phi(t)|^2\,,
  \non
\end{align}
\anold{for a computable constant~$C'$.}
%where we have used the special symbol $C'$ instead of the generic $c$ for a future reference.
\anold{The next estimate follows a similar approach:~from} the second assumption in \eqref{hpFf} and Young's inequality, we obtain~that\anold{, for a suitable constant~$C''$,}
\Beq
  - \nu \iO F(\phi(t)\anold{)}
  \leq \frac 18 \iO |f(\phi(t)|^2
  + C'' \iO |\phi(t)|^2
  + c \,.
  \non
\Eeq
\pier{Now, we observe that the assumptions~\eqref{regphiz} and~\eqref{hppot} allow us to conclude that
\Beq
   \calE \pier{(\Pn \phiz)}
  \leq  c \, , 
  \non
\Eeq
since $\Pn \phiz \to \phiz$ strongly in~$W$ by virtue of~\eqref{convPn}, and uniformly, thanks to the compact embedding $W \hookrightarrow C^0(\overline{\Omega})$.  
\junew{Moreover, we immediately obtain from Young's inequality and~\eqref{regu} that}
\Beq
  \intQt \uu \cdot \vv
  \leq \frac {\lambda_*} 4 \iot \norma{\vv(s)}^2 \, ds
  + c\,.
  \non
\Eeq
We have to} estimate the integral involving~$S$.
Recalling the general notation \eqref{mean}, we write the term under consideration in the form
\Beq
  \intQt \an{h}(\phi) (\mu - \mubar)
  + \intQt \an{h}(\phi) \mubar
  \non
\Eeq
and \juerg{estimate these contributions individually}.
First, by the Poincar\'e and Young inequalities, and using the boundedness assumption on~$\an{h}$, we find that
\Beq
  \intQt \an{h}(\phi) (\mu - \mubar) 
  \leq \CP \iot \norma{\an{h}(\phi(s))} \, \norma{\nabla\mu(s)} \, ds
  \leq \frac \mmin 2 \intQt |\nabla\mu|^2
  + c \,.
  \non
\Eeq
As for the second term,
we recall that the constant functions belong to every~$\Vn$
and test \eqref{quintan} by the constant~$1/|\Omega|$ \anold{to obtain}
\Beq
  \mubar
  = \frac 1 {|\Omega|} \iO \bigl( f'(\phi) + \nu \bigr) \bigl( - \Delta\phi + f(\phi) \bigr) 
  \quad \aet \,.
  \label{mubar}
\Eeq
By also accounting for \pier{\eqref{pier10}} \anold{and the above}, we infer that
\begin{align}
  & \intQt \an{h}(\phi) \mubar
  \leq c \intQt |\mubar|
  = c \iot |\mubar(s)| \, ds
  \non
  \\
  & 
%  \leq  c \intQt \bigl| \bigl( f'(\phi) + \nu \bigr) \bigl( - \Delta\phi + f(\phi) \bigr) \bigr|
 \anold{\leq c }\intQt |{-\Delta\phi }+ f(\phi)|^2
  + c \intQt |f(\phi)|^2 
  + c \,.
  \non
\end{align}
\pier{Finally, it remains to estimate the last term in \eqref{perprimastima}. By exploiting the third assumption in \eqref{hpFf}, we deduce that
\begin{align}
&\pier{{}	- \sigma \intQt \bigl( - \Delta\phi + f(\phi)
\bigr) \bigl(-\Delta \phi + f'(\phi) \phi + \nu \phi \bigr)}
\nonumber
\\
&{}
\leq 
\pier{{} |\sigma|  \intQt \bigl| - \Delta\phi + f(\phi)
\bigr| \bigl(|\Delta \phi| + C_3| f(\phi)| + C_3 + |\nu \phi| \bigr)}
\nonumber
\\
&{}\leq c \intQt |{-\Delta\phi }+ f(\phi)|^2
  + c \intQt |f(\phi)|^2 + c \intQt |\phi|^2
  + c \,.
  \non
\end{align} 
\juerg{Now we come}} back to \eqref{perprimastima}, collect all these inequalities, and rearrange. \juerg{We then} see that the leading terms of the \lhs\ remain, possibly with smaller coefficients,
and that just one expression needs some treatment, namely
\Beq
  \frac 18 \iO |f(\phi(t))|^2
  - (C' + C'') \iO |\phi(t)|^2 \,.
  \non
\Eeq
\pier{However, the growth condition~\eqref{hpf} guarantees the inequality
\begin{equation}
  (C' + C'') s^2 \leq \frac{1}{16} |f(s)|^2 - s^2 + c 
  \qquad \text{for all } s \in \mathbb{R},
  \nonumber
\end{equation}
and, employing it in this form, we are then in a position to apply Gronwall's lemma.}
We conclude, in particular, that
\begin{align}
  & \pier{n^{-1/2} \norma\vvn_{\L\infty\HH}}
  \pier{{}+ \alpha^{1/2} }\norma\vvn_{\L2\VV}
   \pier{{} + \norma{ \eta(\phi_n)^{1/2} | D\vv_n|\/}_{\L2 {\HH^3}}}
  \non
  \\
  & \quad {}
      \pier{{}+ \lambdamin^{1/2} \norma\vvn_{\L2\HH} }
    + \norma{\nabla\mun}_{\anold{\L2 \HH}}
    + \norma{\phin}_{\L\infty H}
    \non
  \\
  & \quad {}
  + \norma{\Delta\phin}_{\L\infty H}
  + \norma{f(\phin)}_{\L\infty H}
  \leq c \,.
  \non
\end{align}
\pier{Then, by using the elliptic regularity estimate \eqref{elliptic1}, we} \juerg{finally arrive at}
\begin{align}
 &\pier{n^{-1/2} \norma\vvn_{\L\infty\HH}}
  \pier{{}+ \alpha^{1/2} \norma\vvn_{\L2\VV}
    + \norma{ \eta(\phi_n)^{1/2} | D\vv_n|\/}_{\L2 {\HH^3}}}
  \nonumber\\
  &\quad{}
  \pier{{}  + \lambdamin^{1/2} \norma\vvn_{\L2\HH} }
  + \norma{\nabla\mun}_{\anold{\L2\HH}}
  + \norma\phin_{\L\infty W}
  \leq c \,.
  \label{primastima}
\end{align}
\pier{Let us point out that the constant $c$ in \eqref{primastima} depends of course on~$\lambdamin$, but is independent of~$\alpha$.}

\step
Consequences

\pier{We recall the computation \eqref{mubar} of the mean value of~$\mu$ and observe that \eqref{primastima}, \eqref{sobolev2} and our regularity assumption \eqref{hppot} on the potential imply that
\Beq
  \norma\phin_\infty + \norma{f(\phin)}_\infty + \norma{f'(\phin)}_\infty
  \leq c \,.
  \label{daprimastima}
\Eeq
Hence, it is straightforward to check that} $\mubar$ is bounded in $\pier{L^\infty (0,T)}$.
By combining this with \eqref{primastima} itself and Poincar\'e's inequality,
we conclude that
\Beq
  \norma\mun_{\L2V} \leq c \,.
  \label{stimamu}
\Eeq

\step
Second a priori estimate

We first observe that \eqref{primastima}, the estimate \eqref{sobolev1} related to the Sobolev embedding $V\emb\Lx6$,
and the \Holder\ inequality imply that
\begin{align}
  & \norma{\vvn\cdot\nabla\phin}_{\L2H}
  \leq c \, \norma{\vvn\cdot\nabla\phin}_{\L2{\Lx3}}
  \non
  \\
  & \leq c \, \norma{\vvn}_{\anold{\L2{\LLx6}}} \, \norma{\nabla\phin}_{\anold{\L\infty{\LLx6}}}
  \leq c \, \norma\vvn_{\L2\VV} \, \norma\phin_{\L\infty W}
  \leq c \,.
  \non
\end{align}
%\comm{[Here, we can also recall that, since $\div \vv =0$, it holds that $\iO \vvn\cdot\nabla\phin z = - \iO \vvn \phin \cdot \nabla z $  and $\vvn \phin \in \L2 \HH$. This of course depends on how we want to formulate the weak form for the phase equation]}
At this point, we take any $\anold{z}\in\L2V$ and test \eqref{secondan} by the projection $\anold{\zn}:=\Pn \anold{z}$.
Then, we integrate over~$(0,T)$ \juerg{to obtain} that
\Beq
  \intQ \dt\phi \, \zn 
  = - \intQ \vv \cdot \nabla\phi \, \zn 
  - \intQ m(\phi) \nabla\mu \cdot \nabla \zn 
  + \intQ
%   \bigl(
  \an{S(\phi)}
%  \pier{{}-\sigma \phi}
%\bigr) 
\zn  \,.
  \non
\Eeq
By virtue of the estimates just obtained,  and owing to Remark~\ref{Projection},
we conclude that
\Beq
  \intQ \dt\phi \, \zn 
  \leq c \, \norma\zn_{\L2V}
  \leq c \, \norma z_{\L2V}\,.
  \non
\Eeq
Since $z\in\L2V$ is arbitrary,  \anold{it is a standard argument to pass to the supremum \pier{with respect to $z$ and infer}} that
\Beq
  \norma{\dt\phin}_{\L2\Vp} \leq c \,.
  \label{secondastima}
\Eeq

\step
Consequence

As announced at the beginning, we have written $T$ instead of the final value $T_n$ 
of the interval where the solution $\soluzn$ to the discrete problem is defined, for simplicity.
However, the above estimates imply that this solution is global.
Indeed, \eqref{secondastima} and the boundedness of $\phin$ in $\L\infty W$ 
with constants that are independent of $T_n$
ensure that $\phin$ is bounded in $\overline\Omega\times[0,T_n)$, and weakly continuous as a function from $[0,T_n]$ to~$W$.
Therefore, if the inequality $T_n<T$ \anold{holds, we could initiate a new Cauchy problem by taking, as the initial value at}
%we could start with a new Cauchy problem 
%by taking as initial value at 
$t=T_n$ the weak limit in $W$ of~$\phin(t)$ as $t$ approaches~$T_n$. \pier{Similar considerations can be repeated for $\vvn$, which is bounded and continuous from $[0,T_n]$ to~$\VV$.}
Since \pier{all} this contradicts the definition of maximal solution, we conclude that $T_n=T$. \anold{Notably, a posteriori, the above estimates are accurate as presented.}

\step
Conclusion

The above estimates and standard weak and weak star compactness results imply that
there exists a triplet $(\vv,\phi,\mu)$ such that\anold{, as $n \to \infty$,}
\begin{align}
   \vvn \to \vv
  & \quad\hbox{weakly in $\pier{\L2\VVz}$}\,,
  \label{convvn}
  \\
   \phin \to \phi
  & \quad\hbox{weakly star in $\H1\Vp\cap\L\infty W$}\,,
  \label{convphin}
  \\
  \mun \to \mu
  & \quad\hbox{weakly in $\L2V$}\,,
  \label{convmun}
\end{align}
at least for a \anold{not} relabeled subsequence.
Moreover, by applying, e.g., \cite[Sect.~8, Cor.~4]{Simon}, and 
\juerg{invoking} the regularity of~$f$,
we deduce the strong convergence \juerg{properties}\pier{%}
\begin{align}
  \phin \to \phi &\quad\hbox{strongly in $\C0{H^{s}(\Omega)}$ for all $s<2$}
  \non\\
  &\quad\hbox{and uniformly in $\overline Q$,}  \label{strongphin}
  \\[2mm]
  \eta(\phin) \to \eta(\phi), &\quad  \lambda (\phin) \to \lambda(\phi), \quad 
   m(\phin) \to m(\phi), \non\\[1mm]
  S(\phin) \to S(\phi), &\quad 
  f(\phin) \to f(\phi),\quad 
  f'(\phin) \to f'(\phi),\non\\
  &\quad \hbox{all uniformly in $\overline Q$, i.e., strongly in $C^0(\overline Q)$}.
  \label{pier2}
\end{align}
We} now prove that the quadruplet $\soluz$, where $w$ is given by~\eqref{quarta}, solves Problem \Pbl\
and satisfies the estimate \eqref{stability}.
\juerg{A part of  the latter is} immediately established:
as a consequence of \eqref{primastima}, \eqref{stimamu}, \eqref{secondastima},
and of the \pier{lower} semicontinuity of norms, we have indeed
\Beq
  \norma\vv_{\L2\VV}
  + \norma\phi_{\H1\Vp\cap\L\infty W}
  + \norma\mu_{\L2V} 
  \leq c \,,
  \label{pier3}
\Eeq
where $c$ has the same dependences as the ones required for $K_1$ in the statement.
Before completing the proof of \eqref{stability}, we have to show that $\soluz$ is the \juerg{sought} solution.

\anold{We defer the treatment of the first equation, which focuses mainly on} $\vv$.

Now, we take any $z\in\L2V$, test \eqref{secondan} by the projection $\zn:=\Pn z$,
and integrate over~$(0,T)$.
By accounting for the strong convergence of $\zn$ to $z$ in $\L2V$, ensured by Remark~\ref{Projection}
(see \eqref{convPnbis} with $Y=V$),
\pier{we let $n$ tend to infinity and observe, in particular, that $\vvn \to \vv$ weakly in $\L2{\HH}$  by \eqref{convvn}\junew{; in addition,} $ (\nabla \phi_n) 
z_n \to  (\nabla \phi) z$ strongly in $\L2{\HH}$ by \eqref{strongphin} and the 
continuity of the embedding $H^{s-1} (\Omega) \emb L^4(\Omega) $ for $s<2$ sufficiently large. Moreover, it turns out that $ \nabla \mu_n \to  \nabla \mu$ weakly in $\L2{\HH}$ by \eqref{convmun} and that
$ m(\phin) \nabla z_n \to
m(\phi) \nabla z$ strongly in $\L2{\HH}$ by \eqref{pier2}. Then, we easily}
conclude that
the time-integrated version of \eqref{seconda} is satisfied with the given test function~$z$.
Since $z$ is arbitrary in $\L2V$, \pier{the time-integrated version of \eqref{seconda}} is equivalent to~\eqref{seconda} \pier{itself.}

\pier{A similar procedure applies to equation \eqref{quintan}, \junew{where this time we take} $z\in\L2W$.
Indeed, we test \eqref{quintan} by the projection $\zn:=\Pn z$, integrate with respect to time and apply \eqref{convPnbis} with $Y=W$ to ensure the strong convergence of $\Delta\zn$ to $\Delta z$ in $\L2H$. \an{T}o pass to the limit as $n\to \infty$, we note\an{, in particular,} that the product $f'(\phin) (-\Delta \phin) $
converges to $f'(\phi) (-\Delta \phi) $ weakly star in $\L\infty H$, due to \eqref{convphin} and \eqref{pier2}. Then, letting $n$ tend to infinity in \eqref{quintan}, we recover the time-integrated version of 
 \begin{align}
 &\iO \bigl( 
    - \Delta\phi + f(\phi)
  \bigr) (-\Delta z)
  + \iO \bigl( f'(\phi) + \nu \bigr) \bigl( - \Delta\phi + f(\phi) \bigr) z
  = \iO \mu z
  \non
  \\
  & \quad \hbox{for every $z\in W$ and \aet}.
  \label{pier4}
\end{align}
Now, if we take $w=-\Delta \phi + f(\phi) $ as a definition of~$w$, \junew{then we infer from~\eqref{pier3}  that}
\Beq
  \norma w_{\L\infty H} \leq c\, .
  \label{pier5}
\Eeq 
At this point, we notice that $w$ solves the variational equality 
\begin{align}
 \iO w (-\Delta z + z ) = 
   \iO \an{g} z
   \quad \hbox{for every $z\in W$ and \aet},
  \label{pier6}
\end{align}
\junew{where the source term $\, \an{g}:= \mu + (1 - f'(\phi) - \nu ) w \, $ is uniformly bounded} in 
$\L2 H$, by \eqref{pier3} and \eqref{pier5}. Then, we can apply~\cite[Lemma~4.2]{CG2} to infer that $w\in \L2 W$, $w$ solves $-\Delta w + w = h $ a.e. in $Q$, and 
\Beq
  \norma w_{\L2 W} \leq c \norma h_{\L2 H}\leq c\, . 
  \label{pier7}
\Eeq 
Hence, we can integrate by parts in the first term of \eqref{pier4} and deduce \eqref{quinta}, 
\junew{at} first for  $z\in W$, and then for all $z\in V$, by density. Obviously,
\eqref{quinta} implies \eqref{terza}: if we write the latter \anold{in} the form
\Beq
	\anold{\iO \nabla w \cdot \nabla z}
%  \iO w (-\Delta z)
  = \iO \bigl( \mu - ( f'(\phi) + \nu ) w \bigr) z
  \quad \hbox{for every $z\in \anold{V}$ and \aet}\,,
  \non
\Eeq
\juerg{then} we recognize the weak formulation of the homogeneous Neumann boundary value problem for 
$-\Delta w =  \mu - ( f'(\phi) + \nu ) w $,
with the forcing term uniformly estimated in $\anold{\L2V}$, as the reader can easily check 
using~\eqref{pier3}, \eqref{hppot} and \eqref{pier7}. 
Therefore, it follows that 
\Beq
  \norma w_{\L2 {\anold{{}\Hx3 \cap W{}}}} \leq c 
  \label{pier7bis}
\Eeq 
by elliptic regularity.}
Hence, coming back to \pier{the definition of~$w$ as $w=-\Delta \phi + f(\phi) $ and to the consequent \junew{identity} \eqref{quarta}, \junew{and} since $f(\phi) $ is uniformly bounded in $\L\infty W$ thanks to \eqref{pier3} and  \eqref{hppot}, \pier{it is not difficult to} deduce that
\Beq
  \norma{\Delta\phi}_{\L2 W} + \norma\phi_{\anold{\L2{\Hx4}}} \leq c \,,
  \label{pier8}
\Eeq
still invoking the elliptic regularity theory. Moreover, now we have that $w-f(\phi) = -\Delta \phi$ is uniformly bounded in in $\L2 {\Hx3}$ \junew{and, consequently, that}  
\Beq
 \norma\phi_{\anold{\L2{\Hx5}}} \leq c \,,
  \label{pier8bis}
\Eeq
which} concludes the proof of \eqref{stability} for the quadruplet $\soluz$  we have found.

\pier{It remains now to verify that $\soluz$} satisfies \eqref{prima} \anold{as well}.
We cannot adopt the argument used for the other equations, 
since we did not prove properties for the projection operator from $\HHz$ onto $\VVzn$
that are analogous to those of~$\Pn$.
Hence, we argue in a different way.
We fix \junew{some $\anold{N}\in\enne$} and assume that $n\geq \anold{N}$.
Then, any $\zz\in\VVzm$ belongs to $\VVzn$ and can be an admissible test function in \eqref{priman}.
Now, take any \juerg{piecewise} linear function $\zz:[0,T]\to\VVzm$ that vanishes at the endpoints.
Then, we can test \eqref{priman} by $\zz$ and integrate over $(0,T)$.
By also integrating by parts in time in the first term, we find that
\Beq
  - \frac 1n \juerg{\intQ} \vvn \cdot \dt\zz 
  + \intQ \bigl(
    \eta(\phin) D\vvn : \nabla\zz
    + \lambda(\phin) \vvn \cdot \zz
  \bigr)
  = \intQ (\mun \nabla\phin + \uu) \cdot \zz\,, 
  \non
\Eeq
and we can let $n$ tend to infinity.
By \eqref{primastima}, we see that the first term tends to zero.
On the other hand, the \pier{remaining} terms tend to those one expects,
and just the product $\mun\nabla\phin$ needs some \anold{further} comment.
\juerg{In view of}~\pier{\eqref{strongphin} and \eqref{convmun},
% and the strong compactness result \cite[Sect.~8, Cor.~4]{Simon} already mentioned,
$\phin$ converges to $\phi$ strongly, e.g., \pier{in $\C0{W^{1,4}(\Omega)}$,}
whence \junew{it follows that} $\mun\nabla\phin$ converges to $\mu\nabla\phi$ weakly in $\L2{\HH}$.}
Therefore, we obtain that
\Beq
  \intQ \bigl(
    \eta(\juerg{\phi}) D\vv : \nabla\zz
    + \lambda(\phi) \juerg{\vv} \cdot \zz
  \bigr)
  = \intQ (\mu \nabla\phi + \uu) \cdot \zz \,,
  \label{perprima}
\Eeq
where $\zz$ is as said.
Since $\anold{N}\in\enne$ was arbitrary and the union of the spaces $\VVzm$ is dense in~$\VVz$,
the set of the piecewise linear functions $\zz$\junew{, which attain their values in this union
and  vanish} at the endpoints of $[0,T]$, 
is dense in $\L2\VVz$.
We conclude that \eqref{perprima} \junew{is valid} for every $\zz\in\L2\VVz$.
But this is equivalent to~\eqref{prima}, 
and the proof of the part of Theorem~\ref{Wellposedness} we wanted to prove is complete.

\Brem
\label{Laplacephi}
\pier{We emphasize that every solution $\soluz$ provided by Theorem~\ref{Wellposedness} enjoys the regularity property (cf.~\eqref{pier8})
\begin{equation}\label{pier9}
  \Delta \phi \in L^2(0,T;W),
\end{equation}
which follows directly from \eqref{quarta}, together with the \junew{regularity properties that}
 $\phi \in L^\infty(0,T;W)$, $w \in L^2(0,T;W)$, and the smoothness of $f$ (see~\eqref{hppot}). This observation will play a  role in the subsequent analysis.}
\Erem

%%%%%%%%%%%%%%%%%%%%%%%%%%%%%%%%%%%%%%%%%%%%%%%%%%%%%%%%%%%%%%%%%%%%%%%%%%%%%%

\section{Proof of Theorem~\ref{Darcy}}
\label{PROOFDARCY}
\setcounter{equation}{0}

%%\begin{proof}
%\noindent\comm{[Andrea: idea of the proof  is that up to the second estimate (eq. (3.24) below, on page 14), it is the same. Namely, we readily have that
%\begin{align*}
%  \norma\vvn_{\L2 \HH}
%  + \norma{\mun}_{\L2 V}
%  + \norma\phin_{\L\infty W}
%  \leq c \,,
%\end{align*}
%with $c$ independent of $n$.
%Next, one notice that
%%\begin{align*}
%%  & \norma{\vvn\cdot\nabla\phin}_{\L2{\Lx {\frac 32}}}
%%  \leq c  \, \norma{\vvn}_{{\L2{\HH}}} \, \norma{\nabla\phin}_{{\L\infty{\LLx6}}}
%% \\ & \quad 
%%  \leq c \, \norma\vvn_{\L2\HH} \, \norma\phin_{\L\infty W}
%%  \leq c \,,
%%\end{align*}
%%and that
% the term in the weak formulation $\iO \vvn \cdot \nabla\phin \, z$ can be rewritten as $ -\iO (\vvn \phin )\cdot  \nabla z$, and that $\norma{\vvn\phin}_{\L2{\HH}} \leq c$, whence also $\dt \phin \in \L2 {\Vp}$.
% Moreover, the strong convergence of $\phin \to \phi$ is enough to let $n\to\infty$ in the nonlinear terms.
% The reconstruction of the chemical potential $w_n$ and the improved regularity for $\phin$ can be 
% carried out as done later on, since they are independent of the regularity of the velocity: it seems 
% that one can even obtain $\phin \in \L2 {\Hx5}$ and $w_n \in \L2 {\Hx3 \cap W}$ uniformly 
% w.r.t.~$n$.
% It is not entirely clear to me how to highlight that, when $n\to\infty$, we have $p=0$ \aeS\, given that we eliminate the pressure. However, perhaps we do not really need anything beyond the weak formulation above.%
%]}%
%%\end{proof}

We start by recalling that, for any $n \in \mathbb{N}$, $(\vvn, \varphi_n, \mu_n, w_n)$ represents a solution \junew{to} Problem~\Pbl\ that corresponds to the function $\eta_n$, \junew{where} 
$$
0<  {\etamin}_{,n} \leq \eta_n(s) \leq {\etamax}_{,n} \quad \hbox{for every $s\in \erre$,}
$$
and the behavior of \junew{the} sequences is such that 
\begin{align}
\|\eta_n\|_{\pier{L^\infty(\mathbb{R})}} \ \hbox{ and, consequently, } \, {\etamin}_{,n} \hbox{ tend to } \, 0 \,  \text{ as } n \to \infty.
\label{etantozero}
\end{align}
The solution $(\vvn, \varphi_n, \mu_n, w_n)$ enjoys the regularity properties \eqref{regv}--\eqref{regw} and satisfies~\Pbl. In particular, \junew{the} conditions \eqref{terza} and \eqref{quarta} can be equivalently expressed in the following strong form:
\begin{align}
 & - \Delta \junew{w_n
  + f'(\phin) w_n
  + \nu w_n 
  = \mun}
  \quad \an{\hbox{in $Q$}\,,}
  \label{Pterza}
  \\
  & -\Delta\junew{\phin
  + f(\phin) 
  = w_n}
   \quad \an{\hbox{in $Q$}\,,}
  \label{Pquarta}  
\end{align} 
which implies that $\varphi_n$ and $\mu_n$ also satisfy the combination of 
\eqref{Pterza}--\eqref{Pquarta}\an{.} Then, in view of \eqref{quinta}), it holds that
\begin{align}
- \Delta \bigl( - \Delta\junew{\phin + f(\phin) \bigr)
+  \bigl( f'(\phin) + \nu \bigr) \bigl( - \Delta\phin + f(\phin) \bigr) = \mun}
 \quad \an{\hbox{in $Q$}\,. }
  \label{Pquinta}
\end{align}
Taking \eqref{regv}--\eqref{regw} into account, together with the regularity assumptions 
on $f$ and $f'$ in \eqref{hppot}, a direct inspection of \eqref{Pquinta} shows that 
all of the terms 
belong to $\L2 V$. Hence, we can reproduce the {\it First a priori estimate} developed in the 
previous section directly at the level of the problem, by testing \eqref{prima} with $\vvn$, 
\eqref{seconda} with $\mun$, and \eqref{Pquinta} with 
$\dt \phin + \sigma \phin \in \L2\Vp$. 

The subsequent computations follow closely those carried out in the previous section and 
lead to the estimate (cf.~\eqref{primastima}--\eqref{stimamu})
\begin{align}
 & \alpha_n^{1/2} \norma\vvn_{\L2\VV}
    + \norma{ \eta_n(\phi_n)^{1/2} | D\vv_n|\/}_{\L2 {\HH^3}}
   + \lambdamin^{1/2} \norma\vvn_{\L2\HH} 
    \nonumber\\
  &\quad{}
  + \norma{\mun}_{\anold{\L2 V}}
  + \norma\phin_{\L\infty W \cap L^\infty (Q)}
  \leq c \,,
  \label{pier11}
\end{align}
where $\alpha_n = \min \graffe{{\etamin}_{,n},\lambdamin}/\CK \to 0 $ as $n\to \infty$ and, importantly, the constant $c$ is uniform with respect to~$n$. 

On the other hand, the {\it Second a priori estimate} of the previous section requires 
a slight modification. Indeed, due to the Sobolev embeddings $V\emb\Lx4$ and $\Lx{4/3} \emb \Vp$,
\junew{we deduce from H\"older's inequality} that 
\begin{align}
  & \norma{\vvn\cdot\nabla\phin}_{\L2{\Vp}}
  \leq c \, \norma{\vvn\cdot\nabla\phin}_{\L2{\Lx{4/3}}}
  \non
  \\
  & \leq c \, \norma{\vvn}_{\anold{\L2{\LLx2}}} \, \norma{\nabla\phin}_{\anold{\L\infty{\LLx4}}}
  \leq c \, \norma\vvn_{\L2\HH} \, \norma\phin_{\L\infty W}
  \leq c \,.
  \non
\end{align}
Hence, by testing \eqref{seconda} with an arbitrary $z \in\L2V$, and using~\eqref{pier11} to control the other terms, we infer that  
\begin{align}
&  \int_0^T  \< \dt\phin, z>\junew{\,dt = \intQ S(\phin) \, z-\intQ m(\phin)\nabla\mun\cdot
\nabla z -\int_0^T\langle \vvn\cdot\nabla\phin,z\rangle\,dt}
\nonumber\\ 
&  \leq c \, \norma z_{\L2V}\,,
  \non
\end{align}
and, consequently, that
\Beq
  \norma{\dt\phin}_{\L2\Vp} \leq c \,.
  \label{pier12}
\Eeq

\junew{Moreover, the estimate \eqref{pier11} 
allows us to deduce that $w_n = -\Delta \phin + f(\phi_n) $ satisfies} 
\Beq
  \norma \wn_{\L\infty H} \leq c\, .
  \label{pier13}
\Eeq 
At this point, we can proceed exactly as in the passage from \eqref{pier6} to \eqref{pier8bis} 
to recover the additional estimate
\Beq 
\norma \wn_{\L2 {\Hx3 \cap W{}}} + \norma{\an{\phi_n}}_{\L2{\Hx5}} \leq c \,, 
\label{pier14} 
\Eeq
which remains valid since these steps are independent of the regularity of $\vvn$. 

Consequently, \junew{from standard weak or weak star compactness arguments, 
we conclude the existence of a subsequence (which is still labeled by $n\in\enne$)
that has  the convergence properties \eqref{convD1}--\eqref{convD5}.} 
In particular, \eqref{convD2} follows from \eqref{pier11} and \eqref{etantozero}. 
Comparing with \junew{the passage-to-the-limit process} carried out in the previous section, 
we see that \eqref{strongphin} and all convergence properties in \eqref{pier2}, except for the first one, 
still hold and allow us to pass to the limit in the present case. 

In particular, \eqref{primaD} and \eqref{secondaD} can be derived by exploiting 
the strong convergence of $\{\phin\}$ together with the density of $\VVz$ in $\HHz$. 
This completes the proof of Theorem~\ref{Darcy}.

%%%%%%%%%%%%%%%%%%%%%%%%%%%%%%%%%%%%%%%%%%%%%%%%%%%%%%%%%%%%%%%%%%%%%%%%%%%%%%

\section{Uniqueness and continuous dependence}
\label{UNIQUENSS}
\setcounter{equation}{0}

This section is devoted to \juerg{the conclusion of}  the proof of Theorem~\ref{Wellposedness}
and to \juerg{the proof of} Theorem~\ref{Contdep}. 
Namely, under the assumption that $\eta$ and $m$ are positive constants,
we prove the uniqueness of the solution to Problem~\Pbl\ 
and the continuous dependence estimate \eqref{contdep}.
\anold{Our approach is as follows: f}irst, we prove an estimate in the direction of \eqref{contdep} 
for any pair of functions $\uu_1$ and $\uu_2$ and arbitrary corresponding solutions
with a constant that also depends on the solutions at hand.
%Nevertheless, we derive uniqueness by considering the particular case $\uu_1=\uu_2$.
%Then, we combine the estimate already proved, the uniqueness of the solution
%and the stability estimate established in the previous section,
%and we conclude the proof of~\eqref{contdep} with a constant $K_2$ 
%that has the properties specified in the statement.
%\juerg{In the following, we} assume that \anold{$\eta\equiv m\equiv 1$} without loss of generality.
\anold{Uniqueness is established by focusing on the particular case $\uu_1 = \uu_2$. We then combine the previously proved estimate, the uniqueness of the solution, and the stability estimate from the preceding section to conclude the proof of~\eqref{contdep}, with a constant $K_2$ possessing the properties stated.  
In what follows, we assume without loss of generality that $\eta \equiv m \equiv 1$.
}

\step
The basic estimate

For given $\uu_1$ and $\uu_2$ satisfying \eqref{regu},
we pick any pair of corresponding solutions $(\vv_1,\phi_1,\mu_1,w_1)$ and $(\vv_2,\phi_2,\mu_2,w_2)$\anold{, and we} set for convenience
\begin{align}
  & \uu := \uu_1 - \uu_2\,, \quad
  \vv := \vv_1 - \vv_2 \,, \quad
  \phi := \phi_1 - \phi_2 \,, \quad
  \mu := \mu_1 - \mu_2 
  \aand
  w := w_1 - w_2 \,\anold{.}
  \non
\end{align}
\anold{We} start proving the estimate mentioned above with a constant that
might depend on the norms of the solutions \anold{involved}.
Indeed, this happens for many of the constants termed~$c$ (or $\cM$, etc.) in the proof,
which also might depend on the \Lip\ constants of $f$, $f'$, $f''$ and $f\!f'$ \anold{encountered during the estimates.}
However, we observe at once that all of the above norms are related to the regularity requirements \Regsoluz\
and that the \Lip\ constants just mentioned depend only on the $L^\infty$ norms 
of the components $\phi_i$ of the solutions we are considering.
This remark is important for the conclusion of our proof.
We introduce two positive parameters $M$ and~$\delta$, whose values are chosen later~on.
\anold{Next}, we take the differences of the equations \eqref{prima}, \eqref{seconda} and \eqref{quinta}, 
written for the data $\uu_1$ and $\uu_2$ and the \anold{corresponding} solutions.
After some rearrangement, we have\anold{,  \aet,}~that
\begin{align}
  & \iO D\vv : \nabla\zz 
\anold{+ \iO \bigl(
   (\lambda(\phi_1) - \lambda(\phi_2)) \vv_1  +\lambda(\phi_2) \vv  \bigr) \cdot \zz
 }%  + \iO \bigl(
%    \lambda(\phi_1) \vv 
%    + (\lambda(\phi_1) - \lambda(\phi_2)) \vv_2 \bigr) \cdot \zz
  \non
  \\
  & = \iO \bigl(
    \mu \nabla\phi_1 + \mu_2 \nabla\phi + \uu
  \bigr) \cdot \zz\,,
  \label{dprima}
\\[2mm]
  & \< \dt\phi , z >
  + \iO \bigl(
    \vv \cdot \nabla\phi_1 + \vv_2 \cdot \nabla\phi
  \bigr) z
  + \iO \nabla\mu \cdot \nabla z
    \non
  \\
   &= \iO \bigl( S(\phi_1) - S(\phi_2) 
%   \pier{{}-\sigma \phi}
\bigr) z\,,
  \label{dseconda}
\\[2mm]
&
\anold{\iO \nabla \bigl( 
    - \Delta\phi + f(\phi_1) - f(\phi_2) 
  \bigr)\cdot \nabla z} 
%\iO \bigl( 
%    - \Delta\phi + f(\phi_1) - f(\phi_2) 
%  \bigr) (-\Delta z)  
  \non
  \\
   &\quad
  + \iO \bigl[ 
    \bigl( f'(\phi_1) - f'(\phi_2) \bigr) (-\Delta\phi_1) 
    + f'(\phi_2) (-\Delta\phi) 
  \bigr] z
  \non
  \\
  & \quad {}
  + \iO \bigl(
    f(\phi_1) f'(\phi_1) - f(\phi_2) f'(\phi_2) 
  \bigr) z 
  \non
  \\
  & \quad {}
    + \nu \iO \bigl(
    -\Delta\phi + f(\phi_1) - f(\phi_2) 
  \bigr) z
  = \iO \mu z\,,
  \label{dquinta}
\end{align}
where $\zz\in \VVz$ is arbitrary in \eqref{dprima}\juerg{,}
$z\in V$ is arbitrary in \eqref{dseconda},  and \juerg{$z\in \anold{V}$ is
arbitrary in  \eqref{dquinta}}.
%At this point, we make a smart choice of the test functions in the above equations.
First, we test \eqref{dprima} \anold{with~}$\vv$.
By recalling the coerciveness inequality \eqref{coercive}, we obtain 
\Beq
  \alpha \, \normaV\vv^2
  \leq \iO |\lambda(\phi_1) - \lambda(\phi_2)| \, \anold{|\vv_1|} \, |\vv|
  + \iO \bigl(
    |\mu| \, |\nabla\phi_1| + |\mu_2| \, |\nabla\phi| + |\uu|
  \bigr) |\vv|\,.
  \label{testdprima}
\Eeq
Next, we test \eqref{dseconda} \anold{with} $M\phi$ to find that
\begin{align}
  & \frac M2 \ddt \, \iO |\phi|^2
  + M \iO \nabla\mu \cdot \nabla \anold{\phi}
  \non
  \\
  & = - M \iO \vv \cdot \nabla\phi_1 \, \phi
  \juerg{\,-\,} M \iO \vv_2 \cdot \nabla\phi \, \phi
  + M \iO \bigl( S(\phi_1) - S(\phi_2)  
%  \pier{{}-\sigma \phi}
  \bigr) \phi \,.
  \label{testdseconda}
\end{align}
\anold{W}e recall the regularity \eqref{regphi} of $\phi_1$ and $\phi_2$, \juerg{and} thus of~$\phi$.
\anold{Therefore, $-M \Delta \phi$ is an admissible test function in \eqref{dseconda}, which leads to}
\begin{align}
  & \frac M2 \ddt \, \iO |\nabla\phi|^2
  + M \iO \nabla\mu \cdot \nabla (-\Delta\phi)
  \non
  \\
  & = - M \iO \vv \cdot \nabla\phi_1 (-\Delta\phi)
  \juerg{\,-\,} M \iO \vv_2 \cdot \nabla\phi (-\Delta\phi)
 \non
  \\
  &\quad{} + M \iO \bigl( S(\phi_1) - S(\phi_2)  
%  \pier{{}-\sigma \phi} 
\bigr) (-\Delta\phi) \,.
  \label{testdsecondaDelta}
\end{align}
\anold{Regarding} \eqref{dquinta}, we observe that the first integral can be written in a different form.
Indeed, we have 
\begin{align}
  & 
  \anold{\iO \nabla\bigl( 
    - \Delta\phi + f(\phi_1) - f(\phi_2) 
  \bigr) \cdot \nabla z}
%\iO \bigl( 
%    - \Delta\phi + f(\phi_1) - f(\phi_2) 
%  \bigr) (-\Delta z)
  = \iO w (-\Delta z)
  \non
  \\
  & = \iO (-\Delta w) z
  = \iO (-\Delta) \bigl( 
    - \Delta\phi + f(\phi_1) - f(\phi_2) 
  \bigr) z\,,
  \non
\end{align}
the integration by parts being allowed since $z\in \anold{V}$ and $w\in\L2W$.
It follows that \eqref{dquinta} can be written as a partial differential equation that holds \aeQ\
and can be multiplied by functions valued in~$H$, in particular, by $-M\Delta\phi+M\Delta^2\phi$ and~$-M\mu$.
By doing this, and integrating over~$\Omega$, we obtain
\begin{align}
  & M \iO |\nabla\Delta\phi|^2
  + M \iO |\Delta^2\phi|^2
  \non
  \\
  & = - M \iO (-\Delta) \bigl( f(\phi_1) - f(\phi_2) \bigr) \bigl( -\Delta\phi + \Delta^2\phi \bigr)
  \non
  \\
  & \quad {}
  - M \iO \bigl( f'(\phi_1) - f'(\phi_2) \bigr) (-\Delta\phi_1) \bigl( -\Delta\phi + \Delta^2\phi \bigr)
  \non
  \\
  & \quad {}
  - M \iO f'(\phi_2) |\Delta\phi|^2 
  - M \iO f'(\phi_2) (-\Delta\phi) \Delta^2\phi
  \non
  \\
  & \quad {}
  - M \iO \bigl(
    f(\phi_1) f'(\phi_1) - f(\phi_2) f'(\phi_2) 
  \bigr) \bigl( -\Delta\phi + \Delta^2\phi \bigr) 
  \non
  \\
  & \quad {}
  - M\nu \iO (-\Delta\phi) \bigl( -\Delta\phi + \Delta^2\phi \bigr)
  - M\nu \iO \bigl( f(\phi_1) - f(\phi_2) \bigr) \bigl( -\Delta\phi + \Delta^2\phi \bigr)
  \non
  \\
  & \quad {}
  + M \iO \nabla\mu \cdot \nabla\phi
  + M \iO \nabla\mu \cdot \nabla(-\Delta\phi)\,, 
  \label{testdquinta}
\end{align}
as well as
\begin{align}
  & M \iO |\mu|^2
  = M \iO (\Delta^2\phi) \mu
  + M \iO (-\Delta) \bigl( 
    f(\phi_1) - f(\phi_2) 
  \bigr) \mu
  \non
  \\
  & \quad {}
  + M \iO \bigl( f'(\phi_1) - f'(\phi_2) \bigr) (-\Delta\phi_1) \mu
  + M \iO f'(\phi_2) (-\Delta\phi) \mu
  \non
  \\
  & \quad {}
  + M \iO \bigl(
    f(\phi_1) f'(\phi_1) - f(\phi_2) f'(\phi_2) 
  \bigr) \mu
  + M \nu \iO \bigl(
    -\Delta\phi + f(\phi_1) - f(\phi_2)  
  \bigr) \mu \,.
  \label{testdquintamu}
\end{align}
At this point, we add relations \accorpa{testdprima}{testdquintamu} to each other and notice \anold{that} some cancellations \anold{occur}.
%Namely, all of the integrals involving $\nabla\mu$ cancel each other.
The \juerg{resulting} \lhs\  is given~by
\begin{align}
  & \alpha \, \normaV\vv^2
  + \frac M2 \ddt \, \iO |\phi|^2
  + \frac M2 \ddt \, \iO |\nabla\phi|^2
  \non
  \\
  & \quad {}
  + M \iO |\nabla\Delta\phi|^2
  + M \iO |\Delta^2\phi|^2
  + M \iO |\mu|^2\,,
  \non
\end{align}
and we have to estimate the terms on the resulting \rhs.
In doing this, we repeatedly make use of some of the inequalities \accorpa{sobolev1}{compact2}, 
as well as of the \Holder\ and Young inequalities \an{(see also  Remark~\ref{Laplacephi})}.
First, we start with the term that is more delicate 
in connection with the choice of the values of the parameters.
It originates from the \rhs\ of \eqref{testdprima} and can be estimated this~way:
\begin{align}
  &  \iO |\mu| \, |\nabla\phi_1| \, |\vv|
  \leq \norma\mu \, \norma{\nabla\phi_1}_4 \, \norma\vv_4
  \leq \CS \normaV\vv \, \norma\mu \, \norma{\nabla\phi_1}_4
  \non
  \\
  & \leq \frac \alpha 8 \, \normaV\vv^2
  + \frac {2\,\CS^2\,R^2} \alpha \, \norma\mu^2 \,,
  \label{critical}
\end{align}
where $\CS$ is the constant appearing in \eqref{sobolev1} 
and $R$ is an upper bound for the norm of $\nabla\phi_1$ in $\anold{\L\infty{\LLx4}}$.
\juerg{In this connection}, we recall that $\phi_1$ belongs to $\L\infty W$ and that $W$ is continuously embedded in~$\Wx{1,4}$.
Next, we have~that
\begin{align}
  & \iO |\lambda(\phi_1) - \lambda(\phi_2)| \, \anold{|\vv_1|} \, |\vv|
 \, \leq \,c \, \norma\phi \, \anold{\norma{\vv_1}_4} \, \norma\vv_4
 \, \leq \,c \, \norma\phi \, \anold{\normaV{\vv_1}} \, \normaV\vv
\non
\\
 & \leq\, \frac \alpha 8 \, \normaV\vv^2 + c \, \anold{\normaV{\vv_1}^2} \, \norma\phi^2\,,
  \non
\end{align}
and we observe at once that the function $t\mapsto\anold{\normaV{\vv_1(t)}^2}$ belongs to $L^1(0,T)$
and is thus suitable for the application of the Gronwall lemma after time integration.
Similar remarks will be omitted in the sequel, for brevity.
For the same reason, we do not recall the regularity of the solutions specified in \Regsoluz\ that
we repeatedly owe to.
For the next term, we have that
\begin{align}
  &  \iO |\mu_2| \, |\nabla\phi| \, |\vv|
  \leq \norma\vv_4 \, \norma{\mu_2}_4 \, \norma{\nabla\phi}
  \leq \CS \, \normaV\vv \, \norma{\mu_2}_4 \, \norma{\nabla\phi}
  \leq \frac \alpha 8 \, \normaV\vv^2
  + c \,  \normaV{\mu_2} ^2\, \norma{\nabla\phi}^2 \,.
  \non
\end{align}
The last term on the \rhs\ of \eqref{testdprima} is trivially estimated by
\Beq
  \iO \uu \cdot \vv
  \leq \frac \alpha 8 \, \normaV\vv^2
  + c \, \norma\uu^2 \,.
  \non
\Eeq
\anold{Concerning} the \rhs\ of \eqref{testdseconda}, we have that
\begin{align}
  & - M \iO \vv \cdot \nabla\phi_1 \, \phi
  \leq M \, \norma\vv_4 \, \norma{\nabla\phi_1}_4 \, \norma\phi
  \leq \frac \alpha 8 \, \normaV\vv^2
  + \cM \, \norma\phi^2\,,
  \non
  \\
  & \juerg{\,-\,}M \iO \vv_2 \cdot \nabla\phi \, \phi
  \leq M \, \norma{\vv_2}_4 \, \norma{\nabla\phi} \, \norma\phi_4
  \leq \cM \, \normaV{\vv_2} \, \normaV\phi^2\,,
  \non
  \\
  & M \iO \bigl( S(\phi_1) - S(\phi_2) 
%  \pier{{}-\sigma \phi} 
\bigr) \phi
  \leq \cM \, \norma\phi^2 \,.
  \non
\end{align}
\anold{Next, we address the \rhs\ of \eqref{testdsecondaDelta}, taking into account} \pier{the 
estimates~\eqref{sobolev1} and~\eqref{compact2}. We deduce that}
\begin{align}
  & - M \iO \vv \cdot \nabla\phi_1 (-\Delta\phi)
  \leq M \norma\vv_4 \, \norma{\nabla\phi_1}_4 \, \norma{\Delta\phi}
  \non
  \\
  & \leq \pier{\frac\alpha 8} \, \normaV\vv^2
  + \pier{\frac {2 \CS^2 M^2}{\alpha}} \, \norma{\nabla\phi_1}_4^2 \, \norma{\Delta\phi}^2
  \non
  \\
  & \leq \pier{\frac\alpha 8} \, \normaV\vv^2
  + \pier{\cM} \, \norma{\Delta\phi}^2
  \leq \pier{\frac\alpha 8} \, \normaV\vv^2
  + \delta \, \norma{\Delta^2\phi}^2
  + \cdeltaM \, \normaV\phi^2 \non
\end{align}
and
\begin{align}
  & \juerg{\,-\,}M \iO \vv_2 \cdot \nabla\phi (-\Delta\phi)
  \leq M \, \norma{\vv_2}_4 \, \norma{\nabla\phi} \, \norma{\Delta\phi}_4
  \leq \cM \, \normaV{\vv_2} \, \norma{\nabla\phi} \, \norma\phi_{\Hx3}
  \non
  \\
  & \leq \cM \, \normaV{\vv_2}^2 \, \norma{\nabla\phi}^2
  + \delta \, \norma{\Delta^2\phi}^2
  + \cdeltaM \, \normaV\phi^2 \,.
  \non
\end{align}
Similarly, \pier{we have that}
\Beq
  M \iO \bigl( S(\phi_1) - S(\phi_2) 
%  \pier{{}-\sigma \phi} 
 \bigr) (-\Delta\phi) 
  \leq c \, M \norma\phi \, \norma{\Delta\phi}
  \leq  \delta \, \norma{\Delta^2\phi}^2
  + \cdeltaM \, \normaV\phi^2 \,.
  \non
\Eeq
\anold{We now proceed to estimate the terms on the \rhs s of \eqref{testdquinta} and \eqref{testdquintamu}. Regarding the first of the latter terms, although straightforward to handle, it is important to note that}
\Beq
  M \iO (\Delta^2\phi) \mu
  \leq \frac M2 \iO |\Delta^2\phi|^2
  + \frac M2 \iO |\mu|^2 \,.
  \non
\Eeq
Next, we consider the first integral on the \rhs\ of \eqref{testdquinta}
and the second one on the \rhs\ of \eqref{testdquintamu}.
Since they can be treated in the same way,
we provide an estimate involving a generic function $z\in\L2H$.
We have
\begin{align}
  & - M \iO (-\Delta) \bigl( f(\phi_1) - f(\phi_2) \bigr) z
  \non
  \\
  & = M \iO \bigl(
    f''(\phi_1) |\nabla\phi_1|^2 - f''(\phi_2) |\nabla\phi_2|^2
    + f'(\phi_1) \Delta\phi_1 - f'(\phi_2) \Delta\phi_2 
  \bigr) z
  \non
  \\
  & = M \iO \bigl[
    \bigl( f''(\phi_1) - f''(\phi_2) \bigr) |\nabla\phi_1|^2
    + f''(\phi_2) \nabla(\phi_1 + \phi_2) \cdot \nabla\phi
  \bigr] z
  \non
  \\
  & \quad {}
  + M \iO \bigl[
    \bigl( f'(\phi_1) - f'(\phi_2) \bigr) \Delta\phi_1
    + f'(\phi_2) \Delta\phi
  \bigr] z
  \non
  \\
  & \leq M c \, \norma\phi_6 \norma{\nabla\phi_1}_6^2 \, \norma z
  + M c \, \norma{\nabla(\phi_1+\phi_2)}_4 \, \norma{\nabla\phi}_4 \, \norma z
  \non
  \\[2mm]
  & \quad {}
  + M c \, \norma\phi_4 \, \norma{\Delta\phi_1}_4 \, \norma z
  + M c \, \norma{\Delta\phi} \, \norma z
  \non
  \\[2mm]
  & \leq \delta \, \norma z^2
  + \cdeltaM \, \norma{\phi_1}_{\Hx2}^2 \, \normaV\phi^2
  + \cdeltaM \, \norma{\phi_1+\phi_2}_{\Hx2}^2 \, \norma\phi_{\Hx2}^2
  \non
  \\[2mm]
  & \quad {}
  + \cdeltaM \, \norma{\phi_1}_{\Hx3}^2 \, \normaV\phi^2
  + \cdeltaM \, \norma{\Delta\phi}^2
  \non
  \\[2mm]
  & \leq \delta \, \norma z^2
  + \cdeltaM \, \norma{\phi_1}_{\Hx3}^2 \, \normaV\phi^2
  + \cdeltaM \, \norma\phi_{\Hx2}^2
  \non
  \\[2mm]
  & \leq \delta \, \norma z^2
  + \cdeltaM \, \norma{\phi_1}_{\Hx3}^2 \, \normaV\phi^2
  + \delta \, \norma{\Delta^2\phi}^2
  + \cdeltaM \, \normaV\phi^2 \,.
  \non
\end{align}
By choosing $z=-\Delta\phi+\Delta^2\phi-\mu$, we derive
an inequality for the sum of the two terms we have to estimate.
Namely, we obtain that
\begin{align}
  & - M \iO (-\Delta) \bigl( f(\phi_1) - f(\phi_2) \bigr) \bigl( -\Delta\phi + \Delta^2\phi \bigr)
  + M \iO (-\Delta) \bigl( f(\phi_1) - f(\phi_2) \bigr) \mu
  \non
  \\
  & \leq \delta \, \norma{-\Delta\phi+\Delta^2\phi-\mu}^2
  + \cdeltaM \, \norma{\phi_1}_{\Hx3}^2 \, \normaV\phi^2
  + \delta \, \norma{\Delta^2\phi}^2
  + \cdeltaM \, \normaV\phi^2 \,.
  \non
\end{align}
On the other hand,
\begin{align}
  & \delta \, \norma{-\Delta\phi+\Delta^2\phi-\mu}^2
  \leq 3\delta \, \norma{\Delta\phi}^2
  + 3\delta \, \norma{\Delta^2\phi}^2
  + 3\delta \, \norma\mu^2
  \non
  \\[2mm]
  & \leq 4\delta \, \norma{\Delta^2\phi}^2
  + \cdelta \, \normaV\phi^2
  + 3\delta \, \norma\mu^2 \,.
  \non
\end{align}
The other integrals on the \rhs\ of \eqref{testdquinta}
that remain after cancellation
can be treated by using the same arguments as before.
Namely, if $I$ is any of them, we have that
\Beq
  I \leq \delta \, \norma{\Delta^2\phi}^2 
  + \psi_1 \, \normaV\phi^2 \,,
  \non 
\Eeq
with some (possibly constant) function $\psi_1\in L^1(0,T)$ depending on $M$ and~$\delta$.
Similarly, if $J$ is any of the other terms \anold{on the \rhs} of \eqref{testdquintamu},
we have that
\Beq
  J \leq \delta \, \norma\mu^2
  + \delta \, \norma{\Delta^2\phi}^2
  + \psi_2 \, \normaV\phi^2 \,,
  \non
\Eeq
\junew{where $\psi_2$ has analogous properties}.
At this \anold{stage}, we can collect all the inequalities we have established
and combine them with the sum of the relations \accorpa{testdprima}{testdquintamu}.
\anold{After straightforward rearrangements,} we \juerg{arrive at} the inequality
\begin{align}
  & \pier{\frac {2\alpha} 8 } \, \normaV\vv^2
  + \frac M2 \ddt \, \normaV\phi^2
  + M \, \norma{\nabla\Delta\phi}^2
  \non
  \\
  & \quad {}
  + \Bigl( \frac M2 - n_1\delta \Bigr) \norma{\Delta^2\phi}^2
  + \Bigl( \frac M2 - \frac {2\,\CS^2\,R^2} \alpha - n_2\delta \Bigr) \norma\mu^2
  \non
  \\
  & \leq \psi \, \normaV\phi^2
  + c \, \norma\uu^2\,,
  \label{lhs}
\end{align}
for some positive integers $n_1$ and $n_2$ and some function $\psi$ belonging to $L^1(0,T)$.
We notice that $\psi$ depends only on the structure of the original system,
upper bounds for some of the norms of the solutions at hand related to \Regsoluz, $M$ and~$\delta$.
Therefore, we can first choose, e.g., $M=8\,\CS^2\,R^2/\alpha$, and then $\delta$ small enough
in order that all the coefficients on the \lhs\ of \eqref{lhs} remain positive.
At this point, we can integrate with respect to time and apply 
 Gronwall's lemma.
\junew{Combining the resulting inequality with elliptic regularity, and invoking  \eqref{elliptic2} and Remark~\ref{Laplacephi}, we finally} conclude that
\Beq
  \norma\vv_{\L2\VV}
  + \norma\phi_{\C0V\cap\L2{\Hx4}}
  + \norma\mu_{\L2H}
  \leq c \, \norma\uu_{\L2\HH}\,\an{.}
  \label{quasicd}
\Eeq
From this, and \eqref{quarta} written for both solutions, we deduce the similar estimate
\Beq
  \norma w_{\L2W}
  \leq c \, \norma\uu_{\L2\HH} \,.
  \label{quasicdw}
\Eeq

\step
Conclusions

Although the constants termed $c$ in \eqref{quasicd} and \eqref{quasicdw} 
depend on the particular solutions we have considered,
these inequalities prove that the solution to Problem \Pbl\ is unique.
Indeed, it suffices to take $\uu_1=\uu_2$ and to recall that the above estimates
hold for arbitrary solutions.
This concludes the proof of Theorem~\ref{Wellposedness}.

\anold{F}or the proof of Theorem~\ref{Contdep}, i.e., of the estimate \eqref{contdep},
we observe that the uniqueness result just established
implies that the solutions involved in \accorpa{quasicd}{quasicdw}
must coincide with the ones we have constructed in the previous section.
Therefore, they satisfy the estimate~\eqref{stability}.
Since the constants that occur in the proof of \accorpa{quasicd}{quasicdw},
in particular those of \accorpa{quasicd}{quasicdw},
depend on the solutions only through upper bounds for some of the norms related to \Regsoluz,
we deduce that they have the same properties as those required for $K_2$ in the statement.

%%%%%%%%%%%%%%%%%%%%%%%%%%%%%%%%%%%%%%%%%%%%%%%%%%%%%%%%%%%%%%%%%%%%%%%%%%%%%%

\section*{Acknowledgements}
\pier{PC was partially supported by the Next Generation EU Project No.~P2022Z7ZAJ 
(``A unitary mathematical framework for modelling muscular dystrophies'').  
AS was partially supported by the ``MUR GRANT Dipartimento di Eccellenza'' 2023--2027 and by the Alexander von Humboldt Foundation.  
Both PC and AS also acknowledge support from the GNAMPA (Gruppo Nazionale per l'Analisi Matematica, la Probabilit\`a e le loro Applicazioni) of INdAM (Isti\-tuto Nazionale di Alta Matematica), with reference project CUP~E53C24001950001.}

%%%%%%%%%%%%%%%%%%%%%%%%%%%%%%%%%%%%%%%%%%%%%%%%%%%%%%%%%%%%%%%%%%%%%%%%%%%%%%

\End{document}

%%%%%%%%%%%%%%%%%%%%%%%%%%%%%%%%%%%%%%%%%%%%%%%%%%%%%%%%%%%%%%%%%%%%%%%%%%%%%%
%%%%%%%%%%%%%%%%%%%%%%%%%%%%%%%%%%%%%%%%%%%%%%%%%%%%%%%%%%%%%%%%%%%%%%%%%%%%%%